\documentclass[reqno,a4paper,oneside]{amsart} 
\usepackage{mathptmx,verbatim,cite,amsmath}
\usepackage{url}

\newtheorem{theorem}{Theorem}[section]
\newtheorem{proposition}[theorem]{Proposition}
\newtheorem{lemma}[theorem]{Lemma}
\newtheorem{corollary}[theorem]{Corollary}

\newtheorem{remark}[theorem]{Remark}
\newtheorem{definition}[theorem]{Definition}

\newcommand{\klg}[1]{\left\{#1\right\}}
\newcommand{\klr}[1]{\left(#1\right)}
\newcommand{\kle}[1]{\left[#1\right]}
\newcommand{\klrk}[1]{(#1)}

\newcommand{\klgk}[1]{\{#1\}}
\newcommand{\betrag}[1]{{\left|#1\right|}}
\newcommand{\betragk}[1]{{|#1|}}
\newcommand{\norm}[1]{\left\|#1\right\|}
\newcommand{\normk}[1]{\|#1\|}

\def\wlim{\mathop{\rm w-}\lim}
\def\toweak{\rightharpoonup}

\def\supp{{\rm supp}}

\def\IR{{\mathbb{R}}}
\def\IN{{\mathbb{N}}}

\DeclareMathOperator{\Cap}{Cap}

\def\domain{\mathsf{D}}
\def\topo{\mathsf{T}}

\def\distance{{\rm{dist}}}

\def\cB{\mathcal{B}}
\def\cD{\mathcal{D}}
\def\cE{\mathcal{E}}
\def\cY{\mathcal{Y}}

\def\cN{\mathcal{N}}
\def\cP{\mathcal{P}}

\def\sC{\;\mathsf{C}}
\def\sE{\;\mathsf{E}}

\def\tsW{\mathcal{W}}
\def\tW{\tilde W}

\def\se{\mathsf{e}}
\def\su{\mathsf{u}}
\def\sg{\mathsf{g}}
\def\sv{\mathsf{v}}
\def\sw{\mathsf{w}}

\begin{document}

\title{The Relative Capacity}
\author{Markus Biegert}
\dedicatory{\upshape%
  Institute of Applied Analysis, University of Ulm, 89069 Ulm, Germany\\
  \texttt{markus.biegert@uni-ulm.de}
}

\begin{abstract}
  The purpose of this article is to introduce the relative $p$-capacity $\Cap_{p,\Omega}$
  with respect to an open set $\Omega$ in $\IR^N$. It is a Choquet capacity on the
  closure of $\Omega$ and extends the classical $p$-capacity $\Cap_p$ in the sense
  that $\Cap_{p,\Omega}=\Cap_p$ if $\Omega=\IR^N$. The importance of the relative
  $p$-capacity stems from the fact that a large class of Sobolev functions defined on a 'bad domain' admits
  a trace on the boundary $\partial\Omega$ which is then unique up to $\Cap_{p,\Omega}$-polar set.
  As an application we prove a characterization of $W^{1,p}_0(\Omega)$ for 
  open sets $\Omega\subset\IR^N$.
\end{abstract}

\maketitle

{\smaller[1]\tableofcontents}


\section{Introduction}

The notion of capacity is fundamental to the analysis of pointwise behavior of Sobolev
functions. Depending on the starting point of the study, the capacity of a set
can be defined in many appropriate ways. The Choquet theory \cite{choquet:54:toc} gives
a standard approach to capacities. Capacity is a necessary tool in classical and
nonlinear potential theory. For example, given an open set $\Omega\subset\IR^N$
the classical $p$-capacity and the relative $p$-capacity can be
used to decide whether a given function $\su\in W^{1,p}(\IR^N)$ lies in
$W^{1,p}_0(\Omega)$ or not. The purpose of this article is to introduce an extension of
the classical $p$-capacity which we call the relative $p$-capacity. Here relative means
with respect to an open and fixed set $\Omega\subset\IR^N$. 

For further results on the classical $p$-capacity and other capacities we refer the
reader to the following books and the references therein:
David R. Adams and Lars I. Hedberg \cite{adams:96:fsp},
Nicolas Bouleau and Francis Hirsch \cite{bouleau:91:df},
Gustave Choquet \cite{choquet:54:toc},
Lawrence C. Evans and Ronald F. Gariepy \cite{evans:92:mtf},
Juha Heinonen and Tero Kilpel{\"a}inen and Olli Martio \cite{heinonen:93:npt},
Jan Mal\'y and William P. Ziemer \cite{ziemer:97:fr} and
Vladimir G. Maz'ya \cite{mazya:85:ssp}.


\section{Preliminaries}\label{sec:rcap}

\subsection{Classical Function Spaces}\label{ss:function}

  Let $\topo$ be a topological space. Then we denote by $C(\topo)$ the space of all
  real-valued and continuous functions on $\topo$ and by $C_c(\topo)$ the subspace
  of $C(\topo)$ consisting of those functions having compact support.
  For an open and non-empty set $\Omega\subset\IR^N$ and $k\in\IN_0$ we let $C^k(\Omega)$
  be the subspace of $C(\Omega)$ consisting of those functions which are $k$ times continuously
  differentiable, that is,
  \[ C^k(\Omega) := \klg{u\in C(\Omega): D^\alpha u\in C(\Omega)\mbox{ for all }\alpha\in\IN_0^N\mbox{ with } \betrag{\alpha}\leq k}.
  \]
  Let $C^\infty(\Omega)$ be the subspace of $C(\Omega)$ given by $C^\infty(\Omega):=\bigcap_{k\in\IN} C^k(\Omega)$
  and let $\cD(\Omega)$ be the space of all {\em test functions} on $\Omega$, that is,
  \[ \cD(\Omega) := C^\infty(\Omega)\cap C_c(\Omega)=\klg{u\in C^\infty(\Omega):\supp(u)\subset\Omega\mbox{ is compact}}
     \subset \cD(\IR^N).
  \]  
  Its topological dual (see Dautray and Lions \cite[Appendix]{dautray:88:man2}) is denoted by $\cD'(\Omega)$ and is called the space of distributions. For
  $p\in[1,\infty)$ the first order Sobolev space $W^{1,p}(\Omega)\subset L^p(\Omega)$ is given by
  \begin{eqnarray*}
     W^{1,p}(\Omega) &:=& \klg{\su\in L^p(\Omega):D^\alpha \su\in L^p(\Omega)
                 \mbox{ in }\cD'(\Omega)\mbox{ for all }\alpha\in\IN_0^N\mbox{ with }\betrag{\alpha}\leq 1} \\
     \norm{\su}^p_{W^{1,p}(\Omega)} &:=& \sum_{\betrag{\alpha}\leq 1} \norm{D^\alpha\su}^p_{L^p(\Omega)}.
  \end{eqnarray*}
  In the following we will work with the closed subspace $\tW^{1,p}(\Omega)$ of the classical Sobolev space
  $W^{1,p}(\Omega)$ defined as the closure of $W^{1,p}(\Omega)\cap C_c(\overline\Omega)$ in $W^{1,p}(\Omega)$
  where the above intersection is defined by
  \[ W^{1,p}(\Omega)\cap C_c(\overline\Omega):=\klg{u|_\Omega:u\in C_c(\overline\Omega), u|_\Omega\in W^{1,p}(\Omega)}.
  \]
  For a real-valued function $u$ we denote by $u^+$ the {\em positive part} and by
  $u^-$ the {\em negative part} of $u$, that is, $u^+:=\max(u,0)=u\vee 0$ and $u^-:=(-u)^+$.

\begin{remark}
  If $\Omega\subset\IR^N$ is an open set with continuous boundary and $p\in[1,\infty)$, then
  the restrictions of functions in $\cD(\IR^N)$ to $\Omega$ are dense in
  $W^{1,p}(\Omega)$ and hence $\tW^{1,p}(\Omega)=W^{1,p}(\Omega)$. See
  Edmunds and Evans \cite[Chap.V, Theorem 4.7]{edmunds:87:std} or
  Maz'ya and Poborchi \cite[Theorem 1.4.2.1]{mazya:97:dfb}.
\end{remark}

\begin{remark}
  For $1<p<\infty$ the space $\tW^{1,p}(\Omega)$ is a uniformly convex (and hence by Milman's theorem a reflexive) Banach space.
  This follows by identifying $\tW^{1,p}(\Omega)$ with a closed subspace of $L^p(\Omega)^{N+1}$.
  For these well-known facts we refer to Alt \cite[Theorem 6.8]{alt:99:lfa},
  Demkowicz and Oden \cite[Proposition 5.13.1(ii)]{oden:96:afa},
  Heuser \cite[Satz 60.4]{heuser:92:fa} and Yosida \cite[Theorem V.2.2]{yosida:80:faa}.
\end{remark}

\subsection{Choquet Capacity}\label{ss:choquet}

 Let $\topo$ be a topological space. For an arbitrary set $\domain$ the power set of $\domain$ is
 denoted by $\cP(\domain)$. A mapping $\sC:\cP(\topo)\to[-\infty,\infty]$ is called a {\em Choquet capacity}
 on $\topo$ if the following properties are satisfied (see Doob \cite[A.II.1]{doob:01:cpt}).
 \begin{itemize}
   \item[(C1)] $\sC$ is increasing; that is, $A\subset B\subset\topo$ implies that $\sC(A)\leq\sC(B)$.
   \item[(C2)] $(A_n)_n\subset\topo$ increasing implies that $\lim_n \sC(A_n)=\sC(\bigcup_n A_n)$.
   \item[(C3)] $(K_n)_n\subset\topo$ decreasing and $K_n$ compact imply $\lim_n \sC(K_n)=\sC(\bigcap_n K_n)$.
 \end{itemize}
 If in addition (C0) holds, then we call $\sC$ a {\em normed Choquet capacity}.
 \begin{itemize}
   \item[(C0)] $\sC(\emptyset)=0$;
 \end{itemize}
 In this case, using (C1), we get that $\sC:\cP(\topo)\to[0,\infty]$.

\subsection{Relative Capacity}\label{ss:rcap}
  Given an open set $\Omega\subset\IR^N$ and $p\in(1,\infty)$ the relative $p$-capacity of an arbitrary set
  $A\subset\overline\Omega$ is defined by
  \[ \Cap_{p,\Omega}(A) := \inf\klg{\norm{\su}^p_{W^{1,p}(\Omega)}:\su\in\cY_{p,\Omega}(A)}
  \]
  where $\cY_{p,\Omega}(A):=\klg{\su\in\tW^{1,p}(\Omega):\exists O\mbox{ open in }\overline\Omega,
          A\subset O,\su\geq 1\mbox{ a.e. on }O\cap\Omega}$.
  Here a.e. is the abbreviation for {\em almost everywhere} with respect to the $N$-dimensional Lebesgue measure.
  In the case $\Omega=\IR^N$ we simply get the classical $p$-capacity which we denote by
  $\Cap_p := \Cap_{p,\IR^N}$.
  The notion of relative $2$-capacity was first introduced by Wolfgang Arendt and
  Mahamadi Warma in \cite{arendt:03:lrb}
  to study the Laplacian with general Robin boundary conditions on arbitrary domains in $\IR^N$.
  Another important application (which will be the subject of a forthcoming paper)
  is the description of vector lattice homomorphisms or isomorphisms between Sobolev spaces.

\section{Properties of the Relative Capacity}
In this section we will systematically collect properties of the relative $p$-capacity. We will
assume throughout the article that $\Omega\subset\IR^N$ is a non-empty open set and $p,q\in (1,\infty)$.

\subsection{Elementary Properties}

\begin{remark}\label{rem:lebesgue}
  It follows directly from the definition that $\lambda^\star(A)\leq \Cap_{p,\Omega}(A)$
  for all sets $A\subset\Omega$ where $\lambda^\star$ denotes the
  outer $N$-dimensional Lebesgue measure.
\end{remark}

\begin{proposition}\label{prop:open-extr}
  Let $O$ be an open set in $\overline\Omega$ (which need not be open in $\IR^N$) and let $p\in(1,\infty)$.
  If $\Cap_{p,\Omega}(O)$ is finite then there is a unique function $\se_O\in\tW^{1,p}(\Omega)$ such that
  \[ \norm{\se_O}_{W^{1,p}(\Omega)}^p=\Cap_{p,\Omega}(O).
  \]
  Moreover, this function satisfies 
  \[ \se_O=1\mbox{ a.e. on } O\cap\Omega\qquad\text{and}\qquad 0\leq \se_O\leq 1\mbox{ a.e. on }\Omega.
  \]
\end{proposition}

\begin{proof}
  The set $\cY_{p,\Omega}(O)=\klg{\su\in\tW^{1,p}(\Omega):\su\geq 1\mbox{ a.e. on }O\cap\Omega}$
  is a closed, convex and non-empty subset of the uniformly convex Banach space $\tW^{1,p}(\Omega)$.
  Let $(\su_n)_n\subset\cY_{p,\Omega}(O)$ be a sequence such that $\norm{\su_n}_{W^{1,p}(\Omega)}^p\to \Cap_{p,\Omega}(O)$.
  By possibly passing to a subsequence we may assume that $\su_n\toweak\se_O$ weakly in $\tW^{1,p}(\Omega)$.
  It follows from Mazur's lemma (see Alt \cite[Lemma 6.13]{alt:99:lfa}) that $\se_O\in\cY_{p,\Omega}(O)$. Moreover, we
  have that $\norm{\se_O}^p_{W^{1,p}(\Omega)}=\Cap_{p,\Omega}(O)$. Using Stampacchia's lemma
  (see Gilbarg and Trudinger \cite[Lemma 7.6 and 7.7]{gilbarg:01:epd}) we get that
  $\se_O':=\se_O^+\wedge 1=\min(\se_O^+,1)\in\cY_{p,\Omega}(O)$ and $\norm{\se_O'}_{W^{1,p}(\Omega)}\leq\norm{\se_O}_{W^{1,p}(\Omega)}$.
  Using that the space $\tW^{1,p}(\Omega)$ is uniformly convex, we get the uniqueness of
  $\se_O$ in $\cY_{p,\Omega}(O)$ and hence $\se_O=\se_O'$ which implies that $\se_O=1$ a.e. on
  $O\cap\Omega$ and $\se_O\in[0,1]$ a.e. on $\Omega$. Note that $\se_O$ is the projection of $0$ onto $\cY_{p,\Omega}(O)$.
\end{proof}

\begin{lemma}\label{lem:compact}
  Let $\topo$ be a Hausdorff space and let $K_n\subset\topo$ be a decreasing sequence of
  compact sets. Then for every open set $V$ containing the intersection $K:=\bigcap_n K_n$
  there exists $n_0\in\IN$ such that $K_n\subset V$ for all $n\geq n_0$.
\end{lemma}

\begin{proof}
  Consider the decreasing sequence of compact sets $(C_n)_n$ given by $C_n:=K_n\setminus V$.
  If $\bigcap_{n=1}^m C_n\not=\emptyset$ for all $m\in\IN$ then $K\setminus V=\bigcap_n C_n\not=\emptyset$
  by Munkres \cite[Theorem 26.9]{munkres:00:top}, a contradiction.
\end{proof}

\begin{theorem}\label{thm:choquet}
  The relative capacity
  $\Cap_{p,\Omega}$ is a normed Choquet capacity on $\overline\Omega$ and for every $A\subset\overline\Omega$
  we have that
  \begin{equation}\label{eq:open}
    \Cap_{p,\Omega}(A) = \inf\klg{\Cap_{p,\Omega}(O):O\mbox{ is open in }\overline\Omega\mbox{ and }A\subset O}.
  \end{equation}
\end{theorem}

\begin{proof}
  That $\Cap_{p,\Omega}$ satisfies the Choquet properties (C0) and (C1) follows immediately
  from the definition
  and the fact that for $A\subset B\subset\overline\Omega$ the inclusion $\cY_{p,\Omega}(B)\subset\cY_{p,\Omega}(A)$ holds.
  The validity of equation \eqref{eq:open} follows also directly from the definition.
  To get the Choquet property (C3) let $(K_n)_n$ be a decreasing sequence of compact subsets of $\overline\Omega$
  and denote by $K$ the intersection of all $K_n$. If $O$ is an open set in $\overline\Omega$ containing $K$ then there
  exists $n_0\in\IN$ such that $K_n\subset O$ for all $n\geq n_0$ (see Lemma \ref{lem:compact}). Hence
  $\Cap_{p,\Omega}(K)\leq\lim_n\Cap_{p,\Omega}(K_n)\leq\Cap_{p,\Omega}(O)$. Taking the infimum over
  all open sets $O$ in $\overline\Omega$ containing $K$ we get by equation \eqref{eq:open} that
  $\Cap_{p,\Omega}(K)=\lim_n \Cap_{p,\Omega}(K_n)$.
  To verify the Choquet property (C2) let $(A_n)_n$ be an increasing sequence of subsets of $\overline\Omega$
  and denote by $A$ the union of all $A_n$. Let $s:=\lim_n \Cap_{p,\Omega}(A_n)\leq\Cap_{p,\Omega}(A)\in[0,\infty]$.
  To get the converse inequality let $\su_n\in\cY_{p,\Omega}(A_n)$ be such that
  $\norm{\su_n}^p_{W^{1,p}(\Omega)}\leq\Cap_{p,\Omega}(A_n)+2^{-n}$.
  We may assume that $s<\infty$, otherwise the equality will be trivial. Therefore $(\su_n)_n$ is a bounded sequence in the
  reflexive Banach space $\tW^{1,p}(\Omega)$ and hence has a weakly convergent subsequence.
  Let $\su\in\tW^{1,p}(\Omega)$ denote the weak limit of this subsequence. By Mazur's lemma there is a
  sequence $(\sv_j)_j$ consisting of convex combinations of the $\su_n$ with $n\geq j$ which
  converges strongly to $\su$. By the triangle inequality we get that
  \[ \norm{\sv_j}^p_{W^{1,p}(\Omega)} \leq \sup_{n\geq j}\norm{\su_n}^p_{W^{1,p}(\Omega)}\leq s+2^{-j}.
  \]
  Moreover, since $\su_n\geq 1$ a.e.\ on $\Omega\cap U_n$ for an open set $U_n$ containing $A_n$ we get that
  there exists an open set $V_n$ (the finite intersection of $U_j$ with $j\geq n$) containing $A_n$
  such that $\sv_n\geq 1$ a.e. on $\Omega\cap V_n$. Since $(\sv_j)_j$ converges to $\su$ we may assume,
  by possibly passing to a subsequence, that $\norm{\sv_{j+1}-\sv_{j}}_{W^{1,p}(\Omega)}\leq 2^{-j}$. Let
  \[ \sw_j:=\sv_j+\sum_{i=j}^\infty \betrag{\sv_{i+1}-\sv_i} \geq \sv_j+\sum_{i=j}^{k-1} (\sv_{i+1}-\sv_i)=\sv_k
     \quad\mbox{ for } k\geq j.
  \]
  Then $\sw_j\in\tW^{1,p}(\Omega)$ and $\sw_j\geq 1$ a.e. $\Omega\cap V$ where the open set $V$ is given
  by $V:=\bigcup_{i=j}^\infty V_i\supset A$. Therefore
  \[ \Cap_{p,\Omega}(A)^{1/p}\leq
         \norm{\sw_j}\leq\norm{\sv_j}+\sum_{i=j}^\infty\norm{\sv_{j+1}-\sv_j}
        \leq (s+2^{-j})^{1/p}+2^{1-j}.
  \]
  For $j\to\infty$ we get that $\Cap_{p,\Omega}(A)\leq s=\lim_n\Cap_{p,\Omega}(A_n)$ which finishes the proof.
\end{proof}

\begin{proposition}\label{prop:cont}
  For a compact set $K\subset\overline\Omega$ we have that
  \begin{eqnarray*}
      \Cap_{p,\Omega}(K) 
        &=& \inf\klg{\norm{u}^p_{W^{1,p}(\Omega)}:u\in \tW^{1,p}(\Omega)\cap C_c(\overline\Omega),u\geq 1\mbox{ on }K}\\
        &=& \inf\klg{\norm{u}^p_{W^{1,p}(\Omega)}:u\in \tW^{1,p}(\Omega)\cap C(\overline\Omega),u\geq 1\mbox{ on }K}.
  \end{eqnarray*}
\end{proposition}

\begin{proof}
  Let $\su\in\cY_{p,\Omega}(K)$ be fixed. Then there exists an open set $U$ in $\IR^N$ containing $K$
  such that $\sv:=(\su\wedge 1)^+=1$ a.e. on $U\cap\Omega$. Let $\eta\in\cD(U)$ be such that $\eta\equiv 1$
  on $K$ and $0\leq\eta\leq 1$ and let $(v_n)_n$ be a sequence in $\tW^{1,p}(\Omega)\cap C_c(\overline\Omega)$
  which converges to $\sv$ in $\tW^{1,p}(\Omega)$. Then $u_n:=\eta+(1-\eta)v_n^+$ converges in
  $\tW^{1,p}(\Omega)$ to $\eta+(1-\eta)\sv=\eta\sv+(1-\eta)\sv=\sv$. Using that
  $u_n\in\tW^{1,p}(\Omega)\cap C_c(\overline\Omega)$, $u_n\geq 1$ on $K$ and
  $\norm{\sv}_{W^{1,p}(\Omega)}\leq\norm{\su}_{W^{1,p}(\Omega)}$ we get that
  \begin{eqnarray*}
     \Cap_{p,\Omega}(K) 
       &\geq& \inf\klg{\norm{u}^p_{W^{1,p}(\Omega)}:u\in\tW^{1,p}(\Omega)\cap C_c(\overline\Omega),u\geq 1\mbox{ on }K}\\
       &\geq& \inf\klg{\norm{u}^p_{W^{1,p}(\Omega)}:u\in\tW^{1,p}(\Omega)\cap C(\overline\Omega),u\geq 1\mbox{ on }K}.
  \end{eqnarray*}
  For the converse inequality we fix a function $u\in\tW^{1,p}(\Omega)\cap C(\overline\Omega)$
  such that $u\geq 1$ on $K$. Then $u_n:=(1+1/n)u\in\cY_{p,\Omega}(K)$ and hence
  \[ \Cap_{p,\Omega}(K) \leq \norm{u_n}^p_{W^{1,p}(\Omega)}\to \norm{u}^p_{W^{1,p}(\Omega)}.
  \]
\end{proof}

\begin{theorem}\label{thm:strongsub}
  The relative $p$-capacity is strongly subadditive, that is, for all $M_1,M_2\subset\overline\Omega$
  \begin{equation}\label{eq:strongsub}
     \Cap_{p,\Omega}(M_1\cup M_2)+\Cap_{p,\Omega}(M_1\cap M_2)\leq\Cap_{p,\Omega}(M_1)+\Cap_{p,\Omega}(M_2)
  \end{equation}
\end{theorem}

\begin{proof}
  Let $\su_j\in\cY_{p,\Omega}(M_j)$ for $j=1,2$ and let $\su:=\max(\su_1,\su_2)$, $\sv:=\min(\su_1,\su_2)$.
  Then $\su\in\cY_{p,\Omega}(M_1\cup M_2)$ and $\sv\in\cY_{p,\Omega}(M_1\cap M_2)$.
  Let $D_1:=\klg{x\in\Omega:\su_1(x)<\su_2(x)}$, $D_2:=\klg{x\in\Omega:\su_1(x)>\su_2(x)}$ and
  $D_3:=\klg{x\in\Omega:\su_1(x)=\su_2(x)}$. By Stampacchia's Lemma
  \begin{eqnarray*}
    \norm{\su}_{W^{1,p}(\Omega)}^p &=& 
        \int_{D_1} \betrag{\su_2}^p+\betrag{\nabla \su_2}^p +
        \int_{D_2} \betrag{\su_1}^p+\betrag{\nabla \su_1}^p +
        \int_{D_3} \betrag{\su_1}^p+\betrag{\nabla \su_1}^p\\
    \norm{\sv}_{W^{1,p}(\Omega)}^p &=& 
        \int_{D_1} \betrag{\su_1}^p+\betrag{\nabla \su_1}^p +
        \int_{D_2} \betrag{\su_2}^p+\betrag{\nabla \su_2}^p +
        \int_{D_3} \betrag{\su_2}^p+\betrag{\nabla \su_2}^p.
  \end{eqnarray*}
  From this we deduce that
  \begin{eqnarray*}
     \Cap_{p,\Omega}(M_1\cup M_2)+\Cap_{p,\Omega}(M_1\cap M_2)
      &\leq& \norm{\su}_{W^{1,p}(\Omega)}^p + \norm{\sv}_{W^{1,p}(\Omega)}^p \\
      &=& \norm{\su_1}_{W^{1,p}(\Omega)}^p + \norm{\su_2}_{W^{1,p}(\Omega)}^p.
  \end{eqnarray*}
  The claim follows now from the definition of the relative $p$-capacity.
\end{proof}

\begin{theorem}
  The relative $p$-capacity is countably subadditive, that is, for all $A_k\subset\overline\Omega$
  \[ \Cap_{p,\Omega}\klr{\bigcup\nolimits_{k\in\IN} A_k} \leq \sum_{k\in\IN} \Cap_{p,\Omega}(A_k).
  \]
\end{theorem}

\begin{proof}
  Let $B_n$ be the union of $A_k$ with $1\leq k\leq n$ and let $A$ be the union of all $A_k$.
  It follows from the strong subadditivity (Theorem \ref{thm:strongsub}) by induction that for all $n\in\IN$
  \[ \Cap_{p,\Omega}\klr{B_n} \leq \sum_{k=1}^n \Cap_{p,\Omega}(A_k).
  \]
  Using the Choquet property (C2) we get
  \[ \Cap_{p,\Omega}\klr{A} = \lim\nolimits_n \Cap_{p,\Omega}\klr{B_n} 
       \leq \lim\nolimits_n \sum_{k=1}^n \Cap_{p,\Omega}(A_k) = \sum_{k\in\IN} \Cap_{p,\Omega}(A_k).
  \]
\end{proof}

\subsection{Relations between Relative Capacities}

\begin{lemma}\label{lem:compacts}
  Let $U\subset V\subset\IR^N$ be non-empty open sets and let $q\leq p$. Then
  for every compact set $K\subset U$ there exists a constant $C=C(U,K,q,p)$ such that
  for all $A\subset K$
  \[ \Cap_{q,V}(A) \leq C\cdot\Cap_{p,U}(A)^{q/p}.
  \]
\end{lemma}

\begin{proof}
  Let $W\subset\subset U$ be an open set which contains $K$ and let $\varphi\in\cD(U)$ be such that
  $\varphi\equiv 1$ on $W$.
  Let $\su\in\cY_{p,U}(A)$ be fixed. Then we define $\sv\in\tW^{1,p}(V)$ by
  $\sv:=\varphi\su$ on $U$ and $\sv:=0$ on $V\setminus U$. Then $\sv\in\cY_{q,V}(A)$ and hence
  by H\"older's inequality
  \[ \Cap_{q,V}(A) \leq \norm{\sv}_{W^{1,q}(V)}^q = \norm{\sv}_{W^{1,q}(U)}^q \leq C_1\norm{\sv}_{W^{1,p}(U)}^q \leq
     C\norm{\su}_{W^{1,p}(U)}^q.
  \]
  Taking the infimum over all $\su\in\cY_{p,U}(A)$ we get the claim.
\end{proof}

\begin{lemma}\label{lem:capesti}
  Let $U\subset V\subset\IR^N$ be non-empty open sets. Then for every
  $A\subset\overline U$ we have that
  \begin{equation}\label{eq:dom}
     \Cap_{p,U}(A) \leq \Cap_{p,V}(A).
  \end{equation}
\end{lemma}

\begin{proof}
  Let $\su\in\cY_{p,V}(A)$. Then $\su|_U\in\cY_{p,U}(A)$ and $\norm{\su|_U}_{W^{1,p}(U)}\leq\norm{\su}_{W^{1,p}(V)}$.
  This implies that $\Cap_{p,U}(A)\leq\Cap_{p,V}(A)$.
\end{proof}

\begin{proposition}\label{prop:polars}
  Let $U\subset V\subset\IR^N$ be non-empty open sets. Then for
  $A\subset U$ we have that
  \begin{equation}\label{eq:polars}
     \Cap_{p,U}(A)=0\quad\Longleftrightarrow\quad \Cap_{p,V}(A)=0.
  \end{equation}
\end{proposition}

\begin{remark}
  In general Equation \eqref{eq:polars} does not hold for $A\subset \partial U\subset\overline U$.
  More precisely, there exists an open, bounded and connected set $\Omega\subset\IR^N$ with $N\geq 2$
  and a smooth set $A\subset\partial\Omega$ such that the $(N-1)$-dimensional Hausdorff measure
  $\mathcal{H}^{N-1}(A)>0$ and $\Cap_{p,\Omega}(A)=0$ for all $p\in(1,\infty)$ - see Biegert \cite[Example 2.5.5]{biegert:05:epvd}.
\end{remark}

\begin{proof}
  From Equation \eqref{eq:dom} we get that $\Cap_{p,V}(A)=0$ implies that $\Cap_{p,U}(A)=0$.
  Hence to prove \eqref{eq:polars} we have to prove the converse implication.
  For this let $\omega_n\subset\subset U$ be an exhausting sequence of $U$ with relatively
  compact sets and assume that $\Cap_{p,U}(A)=0$. Then by Lemma \ref{lem:compacts} there exist
  constants $C_n$ such that
  \[ \Cap_{p,V}(\omega_n\cap A)\leq C_n\cdot \Cap_{p,U}(\omega_n\cap A)=0.
  \]
  Using property (C2) we get that $\Cap_{p,V}(\bigcup_n\omega_n\cap B)=\Cap_{p,V}(B)=0$.
\end{proof}

\begin{definition}
  Let $\Omega\subset\IR^N$ be a domain and let $p\in(1,\infty)$. Then we say that
  $\Omega$ is a {\em $(1,p)$-extension domain} if there exists a bounded linear operator
  $\cE:W^{1,p}(\Omega)\to W^{1,p}(\IR^N)$ such that $\cE\su=\su$ on $\Omega$.
  We say that $\Omega$ has the {\em continuous $(1,p)$-extension property} if there exists
  a bounded linear operator $\cE:W^{1,p}(\Omega)\to W^{1,p}(\IR^N)$ such that
  $\cE\su=\su$ on $\Omega$ and
  $\cE \klr{W^{1,p}(\Omega)\cap C(\overline\Omega)}\subset W^{1,p}(\IR^N)\cap C(\IR^N)$.
\end{definition}

The following is an immediate consequence of Shvartsman \cite{shvartsman:07:eos} and
Haj{\l}asz and Koskela and Tuominen \cite{koskela:08:see}.

\begin{theorem}\label{thm:koskela}
  Let $p\in (1,\infty)$ and $\Omega\subset\IR^N$ be an $(1,p)$-extension domain.
  Then $\Omega$ has the continuous $(1,p)$-extension property.
\end{theorem}

\begin{proof}
  Since $\Omega$ is a $(1,p)$-extension domain, we get from
  \cite[Theorem 2 and Lemma 2.1]{koskela:08:see} that there exists a constant
  $\delta_\Omega>0$ such that $\lambda(B(x,r)\cap\Omega)\geq \delta_\Omega r^N$
  for all $0<r\leq 1$ and $\lambda(\partial\Omega)=0$.
  For a measurable set $A\subset\IR$ we let $M^{1,p}(A)$ be
  the Sobolev-type space introduced by Haj{\l}asz consisting of
  those function $\su\in L^p(A)$ with generalized gradient in $L^p(A)$.
  It follows from \cite[Theorem 1.3]{shvartsman:07:eos} that
  $M^{1,p}(\IR^N)|_{\overline\Omega}=M^{1,p}(\overline\Omega)$ and that
  there exists a linear continuous extension operator
  $\cE:M^{1,p}(\overline\Omega)\to M^{1,p}(\IR^N)$. Using that
  $M^{1,p}(\IR^N)=W^{1,p}(\IR^N)$ as sets with equivalent norms,
  we get that $M^{1,p}(\overline\Omega)=W^{1,p}(\Omega)$
  are equal as sets with equivalent norms, hence the extension
  operator $\cE$ constructed for $M^{1,p}(\overline\Omega)$
  is also a linear continuous extension operator
  from $W^{1,p}(\Omega)$ into $W^{1,p}(\IR^N)$.
   
  To verify that the extension operator $\cE$ constructed by
  Shvartsman maps $W^{1,p}(\Omega)\cap C(\overline\Omega)$
  into $W^{1,p}(\IR^N)\cap C(\IR^N)$ we describe shortly the
  construction of this explicit extension operator, following the arguments
  of Shvartsman.
  There exists a countable family of balls $W=W(\overline\Omega)$ such that
  $\IR^N\setminus\overline\Omega=\bigcup_{B\in W} B$,
  every ball $B=B(x_B,r_B)\in W$ satisfies $3r_B\leq\distance(B,\overline\Omega)\leq 25r_B$ and
  further every point of $\IR^N\setminus\overline\Omega$ is covered by at most $C=C(N)$
  balls from $W$. Let $\klr{\varphi_B}$ for $B\in W$ be a partition of unity associated with
  this Whitney covering $W$ with the properties $0\leq\varphi_B\leq 1$,
  $\supp(\varphi_B)\subset B(x_B,(9/8)r_B)$, $\sum_{B\in W} \varphi_B(x)=1$ on 
  $\IR^N\setminus\overline\Omega$ and $\betrag{\varphi_B(x)-\varphi_B(y)}\leq C\distance(x,y)/r_B$
  for some constant $C>0$ independent of $B$. By \cite[Theorem 2.6]{shvartsman:07:eos} there
  is a family of Borel sets $\klg{H_B:B\in W}$ such that
  $H_B\subset B(x_B,\gamma_1 r_B))\cap\overline\Omega$,
  $\lambda(B)\leq\gamma_2\lambda(H_B)$ for all $B\in W$ with $r_B\leq\delta_\Omega$
  where $\gamma_1$ and $\gamma_2$ are positive constants.
  Now Shvartsman proved \cite[Theorem 1.3 and Equation (1.5)]{shvartsman:07:eos} that
  a continuous extension operator
  $\cE:M^{1,p}(\overline\Omega)\to M^{1,p}(\IR^N)$ is given by
  \[ (\cE\su)(x) := 
       \sum_{B\in W} \su_{H_B}\varphi_B(x)
        \;\text{ for }x\in\IR^N\setminus\overline\Omega,\quad\text{where}\quad
       \su_{H_B}:=\frac{1}{\lambda(H_B)}\int_{H_B} \su\;d\lambda
  \]
 and $(\cE\su)(x):=\su(x)$ if
 $x\in\overline\Omega$. The claim follows now from the construction of $\cE$.
\end{proof}

\begin{theorem}
  Let $U\subset V\subset\IR^N$ be non-empty open sets and assume that $U$ is an $(1,p)$-extension domain.
  Then there exists a constant $C$ depending on $U$ such that for every set $A\subset\overline U$
  \[ \Cap_{p,V}(A)\leq C\cdot \Cap_{p,U}(A) \leq C\cdot \Cap_{p,V}(A).
  \]
\end{theorem}

\begin{proof}
  Let $K\subset\overline U$ be a compact set. By Proposition \ref{prop:cont}
  there exist $u_n\in\tW^{1,p}(U)\cap C_c(\overline U)$ such that
  $u_n\geq 1$ on $K$ and $\norm{u_n}^p_{W^{1,p}(U)}\to \Cap_{p,U}(K)$.
  Let $\cE^\star$ be the extension operator from Theorem \ref{thm:koskela} and define
  $v_n:=\klr{\cE^\star u_n}|_V$. Then $v_n\in\tW^{1,p}(V)\cap C(\overline V)$ and
  $v_n\geq 1$ on $K$. Hence by Proposition \ref{prop:cont} we get that
  \[ \Cap_{p,V}(K)\leq \norm{v_n}_{W^{1,p}(V)}^p \leq 
     \norm{\cE^\star}^p\norm{u_n}^p_{W^{1,p}(U)} \rightarrow
     \norm{\cE^\star}^p\Cap_{p,U}(K).
  \]
  Let $W$ be an open set in $\overline U$. Then there exists an increasing sequence $(K_n)_n$ of
  compact sets such that $\bigcup_n K_n=W$. By the property (C2) we get that
  \[ \Cap_{p,V}(W)=\lim_n \Cap_{p,V}(K_n) \leq 
     \lim_n \norm{\cE^\star}^p\Cap_{p,U}(K_n)=\norm{\cE^\star}^p\Cap_{p,U}(W).
  \]
  Now let $A\subset\overline U$ be arbitrary. Then by Theorem \ref{thm:choquet}
  \begin{eqnarray*}
    \Cap_{p,V}(A) 
        &=& \inf\klg{\Cap_{p,V}(O):O\mbox{ is open in }\overline V\mbox{ and } A\subset O} \\
        &=& \inf\klg{\Cap_{p,V}(O\cap\overline U):O\mbox{ is open in }\overline V\mbox{ and } A\subset O} \\
        &=& \inf\klg{\Cap_{p,V}(W):W\mbox{ is open in }\overline U\mbox{ and } A\subset W} \\
        &\leq& \norm{\cE^\star}^p \inf\klg{\Cap_{p,U}(W):W\mbox{ is open in }\overline U\mbox{ and }
               A\subset W}\\
        &=& \norm{\cE^\star}^p\Cap_{p,U}(A).
  \end{eqnarray*}
  The remaining inequality follows from Lemma \ref{lem:capesti}.
\end{proof}

\begin{theorem}
  For $q\leq p$ and $A\subset\overline\Omega$ with $\Cap_{p,\Omega}(A)=0$ we have that $\Cap_{q,\Omega}(A)=0$.
\end{theorem}

\begin{proof}
  Let $n\in\IN$, $\su\in\cY_{p,\Omega}(A\cap B(0,n))$ and $\eta\in C^1_c(B(0,2n))$ be such that
  $\eta\equiv 1$ on $B(0,n)$. Then we get by H\"older's inequality
  \[ \Cap_{q,\Omega}(A\cap B(0,n))\leq \norm{\su\eta}^q_{W^{1,q}(\Omega)} \leq
      C_1\cdot\norm{\su\eta}^q_{W^{1,p}(\Omega)} \leq C\cdot\norm{\su}^q_{W^{1,p}(\Omega)}
  \]
  where $C$ is a constant independent of $\su$. Taking the infimum over all such $\su$ we get that
  \[ \Cap_{q,\Omega}(A\cap B(0,n))\leq C\cdot\Cap_{p,\Omega}(A\cap B(0,n))^{q/p} = 0
  \]
  and hence by the property (C2) we get that $\Cap_{q,\Omega}(A)=0$.
\end{proof}

\subsection{Quasicontinuity and Polar Sets}

The aim of this subsection is to prove the existence and uniqueness of $\Cap_{p,\Omega}$-quasi continuous
representatives on $\overline\Omega$ for Sobolev functions in $\tW^{1,p}(\Omega)$.

\begin{definition}
  A set $P\subset\overline\Omega$ is said to be $\Cap_{p,\Omega}$-polar if
  $\Cap_{p,\Omega}(P)=0$. A pointwise defined function $u$ on $\domain\subset\overline\Omega$
  is called $\Cap_{p,\Omega}$-quasi continuous on $\domain$ if for each $\varepsilon>0$
  there exists an open set $V$ in $\overline\Omega$ with $\Cap_{p,\Omega}(V)<\varepsilon$ such that
  $u$ restricted to $\domain\setminus V$ is continuous. We say that a property
  holds $\Cap_{p,\Omega}$-quasi everywhere (briefly $\Cap_{p,\Omega}$-q.e.) if
  it holds except for a $\Cap_{p,\Omega}$-polar set.
\end{definition}

\begin{lemma}\label{lem:conver}
  If $\su\in\tW^{1,p}(\Omega)$ and $u_k\in\tW^{1,p}(\Omega)\cap C_c(\overline\Omega)$
  is such that
  \[ \sum_{k=1}^\infty 2^{kp}\norm{\su-u_k}^p_{W^{1,p}(\Omega)}<\infty,
  \]
  then the pointwise limit $\tilde u:=\lim_k u_k$ exists $\Cap_{p,\Omega}$-quasi everywhere
  on $\overline\Omega$, $\tilde u:\overline\Omega\to\IR$ is $\Cap_{p,\Omega}$-quasi continuous 
  and $\tilde u=\su$ almost everywhere on $\Omega$.
\end{lemma}

\begin{proof}
  Let $G_k:=\klgk{x\in\overline\Omega:\betragk{u_{k+1}(x)-u_k(x)}>2^{-k}}$. Then $G_k$ is an open set
  in $\overline\Omega$ and $2^k\betragk{u_{k+1}-u_k}\geq 1$ on $G_k$. It follows that
  \[ \Cap_{p,\Omega}(G_k) \leq 2^{kp}\norm{u_{k+1}-u_k}^p_{W^{1,p}(\Omega)}
  \]
  and hence $\sum_{k} \Cap_{p,\Omega}(G_k)<\infty$. Given $\varepsilon>0$ there exists $k_0\in\IN$
  such that $\Cap_{p,\Omega}(G)<\varepsilon$ where $G:=\bigcup_{k\geq k_0} G_k$.
  Since $\betrag{u_{k+1}-u_k}\leq 2^{-k}$ on $\overline\Omega\setminus G$ for all $k\geq k_0$ we have
  that $(u_k)_k$ is a sequence of continuous functions on $\overline\Omega$ which converges
  uniformly on $\overline\Omega\setminus G$. Since $\varepsilon>0$ was arbitrary
  we get that $\tilde u:=\lim_k u_k$ exists $\Cap_{p,\Omega}$-quasi everywhere on $\overline\Omega$
  and $\tilde u|_{\overline\Omega\setminus G}$ is continuous. To see that $\tilde u$ coincides with
  $\su$ almost everywhere on $\Omega$ we argue as follows. Since $u_k$ converges to $\su$
  in $W^{1,p}(\Omega)$ (by possibly passing to a subsequence) we have that $u_k$ converges to
  $\su$ almost everywhere. Since $(u_k)_k$ converges to $\tilde u$ $\Cap_{p,\Omega}$-quasi everywhere
  on $\overline\Omega$ (and hence almost everywhere on $\Omega$) we get that
  $\tilde u=\su$ almost everywhere on $\Omega$ (see Remark \ref{rem:lebesgue}).
\end{proof}

\begin{theorem}\label{thm:cont-version}
  For every $\su\in\tW^{1,p}(\Omega)$ there exists a $\Cap_{p,\Omega}$-quasi
  continuous function $\tilde u:\overline\Omega\to\IR$ such that
  $\tilde u=\su$ almost everywhere on $\Omega$, that is, $\tilde u\in\su$.
\end{theorem}

\begin{proof}
  Let $\su\in\tW^{1,p}(\Omega)$. Then by definition there exists a sequence
  $u_n\in\tW^{1,p}(\Omega)\cap C_c(\overline\Omega)$ such that $u_n\to \su$ in
  $\tW^{1,p}(\Omega)$. By possibly passing to a subsequence the sequence
  $(u_n)_n$ satisfies the assumptions of Lemma \ref{lem:conver}.
\end{proof}

\begin{lemma}\label{lem:approx}
  Let $A\subset\overline\Omega$ and let $u\in\su\in\tW^{1,p}(\Omega)$ be a
  $\Cap_{p,\Omega}$-quasi continuous version of $\su$ such that
  $u\geq 1$ $\Cap_{p,\Omega}$-quasi everywhere on $A$. Then there
  is a sequence $(\su_n)_n\subset \cY_{p,\Omega}(A)$ which converges to $\su$
  in $\tW^{1,p}(\Omega)$.
\end{lemma}

\begin{proof}
  Let $\varepsilon>0$ and let $k\in\IN$ be such that $\norm{\su-w_k}_{W^{1,p}(\Omega)}\leq\varepsilon$
  where $w_k$ is given by $w_k:=\max(u,-k)$. Let $V$ be an open set in $\overline\Omega$
  such that $w_k$ restricted to $\overline\Omega\setminus V$ is continuous,
  $\sw_k\geq 1$ everywhere on $A\setminus V$ and $\Cap_{p,\Omega}(V)\leq (k+1+\varepsilon k)^{-p}\varepsilon^p$.
  Let $\psi$ be a capacitary extremal for $V$ (see Proposition \ref{prop:open-extr}) 
  and let $\sv_\varepsilon:=(1+\varepsilon)w_k+(k+1+\varepsilon k)\psi\geq\psi$. Then
  \begin{eqnarray*}
     \norm{\su-\sv_\varepsilon}_{W^{1,p}(\Omega)} 
   &\leq& \varepsilon+\norm{w_k-\sv_\varepsilon}_{W^{1,p}(\Omega)} \leq
          \varepsilon+\varepsilon\norm{w_k}_{W^{1,p}(\Omega)}+\varepsilon \\
   &\leq& \varepsilon\klrk{2+\norm{\su}_{W^{1,p}(\Omega)}+\varepsilon}.
  \end{eqnarray*}
  For the open set $G:=V\cup \klg{x\in\overline\Omega\setminus V:w_k(x)>1/(1+\varepsilon)}$ in $\overline\Omega$
  we have that $\sv_\varepsilon\geq 1$ a.e. on $\Omega\cap G$ and $A\subset G$. Hence
  $\sv_\varepsilon\in\cY_{p,\Omega}(A)$. The claim follows with $\su_n:=\sv_{1/n}$.
\end{proof}

\begin{lemma}\label{lem:normest}
  Let $u\in\su\in\tW^{1,p}(\Omega)$ be a $\Cap_{p,\Omega}$-quasi continuous version of $\su$ and let
  $a\in(0,\infty)$. Then
  \[ \Cap_{p,\Omega}\klrk{\klgk{x\in\overline\Omega:u(x)>a}}\leq a^{-p}\norm{\su^+}^p_{W^{1,p}(\Omega)}.
  \]
\end{lemma}

\begin{proof}
  Let $A:=\klgk{x\in\overline\Omega:u(x)>a}$. By Lemma \ref{lem:approx} there exists a sequence
  $(\su_n)_n\in\cY_{p,\Omega}(A)$ which converges to $a^{-1}\su^+$ in $\tW^{1,p}(\Omega)$.
  Note that $u^+$ is a $\Cap_{p,\Omega}$-quasi continuous version of $\su^+$. Hence
  \[ \Cap_{p,\Omega}\klrk{\klg{x\in\overline\Omega:u(x)>a}}\leq\norm{\su_n}^p_{W^{1,p}(\Omega)}\to
     a^{-p}\norm{\su^+}^p_{W^{1,p}(\Omega)}.
  \]
\end{proof}

\begin{theorem}\label{thm:esti}
  Let $\su,\sv\in\tW^{1,p}(\Omega)$ be such that $\su\leq\sv$ a.e. on $U\cap\Omega$ where $U$ is an open set in $\overline\Omega$.
  If $u\in\su$ and $v\in\sv$ are $\Cap_{p,\Omega}$-quasi continuous versions of $\su$ and $\sv$, respectively,
  then $u\leq v$ $\Cap_{p,\Omega}$-quasi everywhere on $U$.
\end{theorem}

\begin{proof}
  Let $W$ be an open set in $\IR^N$ such that $U=W\cap\overline\Omega$ and let $(K_n)_n$ be a
  sequence of compact sets such that $U=\bigcup_n K_n$. For the sequence of compact sets
  we choose $\varphi_n\in\cD(W)$ non-negative such that $\varphi_n\equiv 1$ on $K_n$.
  Then the function $w_n:=\varphi_n(u-v)^+=0$ a.e. on $\Omega$ and we get by
  Lemma \ref{lem:normest}, using that $\varphi_n (u-v)^+$ is $\Cap_{p,\Omega}$-quasi
  continuous, that $w_n=0$ $\Cap_{p,\Omega}$-quasi everywhere on $\overline\Omega$
  and hence that $u\leq v$ $\Cap_{p,\Omega}$-quasi everywhere on $K_n$ for each $n\in\IN$.
  Since the countable union of $\Cap_{p,\Omega}$-polar sets is $\Cap_{p,\Omega}$-polar
  we get that $u\leq v$ $\Cap_{p,\Omega}$-quasi everywhere on $U$.
\end{proof}

\begin{theorem}\label{thm:quasicontrep}
  Let $\su\in\tW^{1,p}(\Omega)$. Then there exists a unique (up to a $\Cap_{p,\Omega}$-polar set)
  $\Cap_{p,\Omega}$-quasi continuous function $\tilde u:\overline\Omega\to\IR$ such that
  $\su=\tilde u$ a.e. on $\Omega$.
\end{theorem}

\begin{proof}
  The existence follows from Theorem \ref{thm:cont-version}. To show uniqueness we let
  $u_1,u_2\in\su\in\tW^{1,p}(\Omega)$ be two quasi-continuous versions $\su$.
  Then $u_1=u_2$ a.e. on $\Omega$ and hence by Theorem \ref{thm:esti} we get that
  $u_1=u_2$ $\Cap_{p,\Omega}$-quasi everywhere on $\overline\Omega$.
\end{proof}

\begin{definition}
  By $\cN_p^\star(\Omega)$ we denote the set of all $\Cap_{p,\Omega}$-polar sets in $\overline\Omega$
  and we denote by $C_p(\Omega)$ the space of all $\Cap_{p,\Omega}$-quasi continuous functions
  $u:\overline\Omega\to\IR$. On $C_p(\Omega)$ we define the equivalence relation $\sim$ by
  \[ u\sim v\quad :\Leftrightarrow \quad \exists P\in\cN_p^\star(\Omega): u=v \mbox{ everywhere on } \overline\Omega\setminus P.
  \]
  For a function $u\in C_p(\Omega)$ we denote by $[u]$ the equivalence class of $u$ with respect to $\sim$.
  Now the refined Sobolev space $\tsW^{1,p}(\Omega)$ is defined by
  \[ \tsW^{1,p}(\Omega) := \klg{\tilde\su:\su\in\tW^{1,p}(\Omega)}\subset C_p(\Omega)/\sim
  \]
  where $\tilde\su:=[u]$ with $u\in\su\in\tW^{1,p}(\Omega)$ $\Cap_{p,\Omega}$-quasi continuous.
  We equip $\tsW^{1,p}(\Omega)$ with the norm $\norm{\cdot}_{W^{1,p}(\Omega)}$. Note
  that it is isometrically isomorphic to $\tW^{1,p}(\Omega)$ by Theorem \ref{thm:quasicontrep}.

  For a sequence $(\su_n)_n$ in $\tsW^{1,p}(\Omega)$ and $\su\in\tsW^{1,p}(\Omega)$
  we say that {\em $(\su_n)_n$ converges $\Cap_{p,\Omega}$-quasi everywhere to $\su$} if
  for every $u_n\in\su_n$ and $u\in\su$ there exists a $\Cap_{p,\Omega}$-polar set $P$ such that
  $u_n\to u$ everywhere on $\overline\Omega\setminus P$.
  We say that {\em $(\su_n)_n$ converges $\Cap_{p,\Omega}$-quasi uniformly to $\su$}
  if for every $u_n\in\su_n$, $u\in\su$ and $\varepsilon>0$ there exists an open set $G$ in $\overline\Omega$
  such that $\Cap_{p,\Omega}(G)\leq\varepsilon$ and $u_k\to u$ uniformly (everywhere) on
  $\overline\Omega\setminus G$.
\end{definition}

\begin{theorem}\label{thm:convsubse}
  If $\su_n\in\tsW^{1,p}(\Omega)$ converges to $\su\in\tsW^{1,p}(\Omega)$ in
  $\tsW^{1,p}(\Omega)$, then there exists a subsequence which converges $\Cap_{p,\Omega}$-quasi
  everywhere and -quasi uniformly on $\overline\Omega$ to $\su$.
\end{theorem}

\begin{proof}
  Let $(\su_{n_k})_k$ be a subsequence of $(\su_n)_n$ such that
  $\sum_{k\in\IN} k^p\normk{\su_{n_k}-\su}^p_{\tsW^{1,p}(\Omega)}<\infty$.
  We show that this subsequence converges $\Cap_{p,\Omega}$-quasi everywhere 
  and -quasi uniformly on $\overline\Omega$ to $\su$.
  Let $u_{n_k}\in\su_{n_k}$ and $u\in\su$ be fixed and define
  $G_k:=\klg{x\in\overline\Omega:\betrag{u_{n_k}(x)-u(x)}>k^{-1}}$.
  We show that $u_{n_k}(x)\to u(x)$ for all $x\in\overline\Omega\setminus P$
  where $P:=\bigcap_{j=1}^\infty \bigcup_{k=j}^\infty G_k$.
  If $x\in\overline\Omega\setminus P$ then there exists $j_0\in\IN$ such that
  $x\not\in \bigcup_{k=j_0}^\infty G_k$, that is, $\betrag{u_{n_k}(x)-u(x)}\leq k^{-1}$
  for all $k\geq j_0$ and hence $u_{n_k}(x)\to u(x)$
  uniformly on $\overline\Omega\setminus\bigcup_{k=j_0}^\infty G_k$ and
  everywhere on $\overline\Omega\setminus P$. We show that $P$ is a $\Cap_{p,\Omega}$-polar set.
  Let $\varepsilon>0$ be given. Then there exists $N=N(\varepsilon)$ such that
  $\sum_{k=N}^\infty k^p\norm{\su_{n_k}-\su}\leq\varepsilon$. By Lemma \ref{lem:normest} we get that
  $\Cap_{p,\Omega}(G_k)\leq k^p\norm{\su_{n_k}-\su}_{\tsW^{1,p}(\Omega)}^p$ and hence
  $\Cap_{p,\Omega}(\bigcup_{k=N}^\infty G_k)\leq\varepsilon$. Therefore $\Cap_{p,\Omega}(P)\leq\varepsilon$
  and since $\varepsilon>0$ was arbitrary the claim follows.
\end{proof}

\begin{lemma}\label{lem:loccont}
  Let $U\subset\Omega\subset\IR^N$ be non-empty open sets. A function $u:\Omega\to\IR$ is
  $\Cap_{p,\Omega}$-quasi continuous on $U$ if and only if
  $u$ is $\Cap_{p,\Omega}$-quasi continuous on every set $\omega\subset\subset U$.
\end{lemma}

\begin{proof}
  Assume that $u$ is $\Cap_{p,\Omega}$-quasi continuous on every set $\omega\subset\subset U$.
  Let $\omega_n\subset\subset U$ be such that $\bigcup_n\omega_n=U$ and let $\varepsilon>0$ be given.
  Then there exists an open set $V_n\subset\omega_n$ such that $u|_{\omega_n\setminus V_n}$ is continuous
  and $\Cap_{p,\Omega}(V_n)\leq \varepsilon 2^{-n}$. Let $V:=\bigcup_n V_n$. Then
  $\Cap_{p,\Omega}(V)\leq \sum_n \Cap_{p,\Omega}(V_n)\leq\varepsilon$ and $u|_{U\setminus V}$ is continuous.
  In fact, if $x\in U\setminus V$ then there exists $n_0\in\IN$ such that $x\in\omega_{n_0}\setminus V_{n_0}$.
  If $(x_k)_k$ is a sequence in $U\setminus V$ converging to $x$ then there exists $k_0$ such
  that $x_k\in\omega_{n_0}\setminus V_{n_0}$ for all $k\geq k_0$. Since $u|_{\omega_{n_0}\setminus V_{n_0}}$ is continuous
  we get that $u(x_k)\to u(x)$ as $k\to\infty$ and hence that $u|_{U\setminus V}$ is continuous.
\end{proof}

\begin{theorem}\label{thm:quasicont}
  Let $U\subset V\subset\IR^N$ be non-empty open sets and let
  $u$ be a function from $U$ into $\IR$. Then
  $u$ is $\Cap_{p,U}$-quasi continuous if and only if
  $u$ is $\Cap_{p,V}$-quasi continuous.
\end{theorem}

\begin{proof}
  If $u$ is $\Cap_{p,V}$-quasi continuous, then $u$ is $\Cap_{p,U}$-quasi continuous
  by Lemma \ref{lem:capesti}. Assume now that $u$ is $\Cap_{p,U}$-quasi continuous and let
  $\omega\subset\subset U$ be a relatively compact set in $U$. We show that $u$
  is $\Cap_{p,V}$-quasi continuous on $\omega$. For this let $\varepsilon>0$ be fixed.
  By Lemma \ref{lem:compacts} there exists a constant $C>0$ such that
  $\Cap_{p,V}(A)\leq C\cdot\Cap_{p,U}(A)$ for all $A\subset\omega$. Since
  $u$ is $\Cap_{p,U}$-quasi continuous on $\omega$ there exists an open set
  $W\subset\omega$ with $\Cap_{p,U}(W)\leq\varepsilon/C$ such that $u|_{\omega\setminus W}$ is
  continuous. Since $\Cap_{p,V}(W)\leq C\Cap_{p,U}(W)\leq\varepsilon$ and $\varepsilon>0$
  was arbitrary we get that $u$ is $\Cap_{p,V}$-quasi continuous on $\omega$.
  Since $\omega\subset\subset U$ was arbitrary we get by Lemma \ref{lem:loccont} that 
  $u$ is $\Cap_{p,V}$-quasi continuous.
\end{proof}

\begin{corollary}
  Let $u\in\su\in\tW^{1,p}(\Omega)$ be a $\Cap_p$-quasi continuous function.
  Then $u=\tilde\su$ $\Cap_{p}$-quasi everywhere on $\Omega$. 
\end{corollary}

\begin{proof}
  Let $v\in\tilde\su$. By Theorem \ref{thm:quasicont} $v$ is $\Cap_p$-quasi continuous
  on $\Omega$. Let $\omega_n\subset\subset\Omega$ be an increasing sequence of
  relatively compact sets in $\Omega$ such that $\bigcup_n \omega_n=\Omega$.
  Let $\varphi_n\in\cD(\Omega)$ be such that $\varphi_n\equiv 1$ on $\omega_n$.
  Since $v=u$ a.e. on $\Omega$ we get that $\varphi_n v=\varphi_n u$ a.e. on
  $\Omega$. Since $\varphi_n v,\varphi_n u\in W^{1,p}(\IR^N)$ are $\Cap_p$-quasi
  continuous on $\IR^N$ we get by Theorem \ref{thm:esti} that
  $\varphi_n v=\varphi_n u$ $\Cap_p$-q.e. on $\IR^N$ and hence
  $v=u$ $\Cap_p$-q.e. on $\omega_n$. Since $(\omega_n)_n$ was exhausting
  we get that $v=u$ $\Cap_p$-quasi everywhere on $\Omega$.
\end{proof}

\subsection{Capacitary Extremals}

The aim of this subsection is to prove the existence and uniqueness of capacitary extremals
and to characterize them.

\begin{theorem}\label{thm:ypo}
  Let $A\subset\overline\Omega$ and $\su\in\tW^{1,p}(\Omega)$. Then $\su\in\overline\cY_{p,\Omega}(A)$
  if and only if $\tilde\su\geq 1$ $\Cap_{p,\Omega}$-q.e. on $A$.
\end{theorem}

\begin{proof}
  If $\tilde\su\geq 1$ $\Cap_{p,\Omega}$-q.e. on $A$ then $\su\in\overline\cY_{p,\Omega}(A)$ by
  Lemma \ref{lem:approx}. For the converse implication let $\su\in\overline\cY_{p,\Omega}(A)$.
  By Theorem \ref{thm:convsubse} there exists a sequence
  $(\su_n)_n\subset\cY_{p,\Omega}(A)$ such that $\tilde\su_n\to\tilde\su$
  $\Cap_{p,\Omega}$-q.e. on $\overline\Omega$. For every $n\in\IN$ there exists
  an open set $O_n$ in $\overline\Omega$ containing $A$ such that
  $\su_n\geq 1$ a.e. on $\Omega\cap O_n$. Hence $\tilde\su_n\geq 1$ $\Cap_{p,\Omega}$-q.e.
  on $A$ by Theorem \ref{thm:esti}. This shows that $\tilde\su\geq 1$ $\Cap_{p,\Omega}$-q.e.
  on $A$.
\end{proof}

\begin{theorem}\label{thm:quasi-cap}
  For $A\subset\overline\Omega$ the relative $p$-capacity of $A$ is given by
  \begin{eqnarray}
     \Cap_{p,\Omega}(A) 
     &=& \inf\klg{\norm{\su}^p_{W^{1,p}(\Omega)}:\su\in\tsW^{1,p}(\Omega),\su\geq 1\mbox{ $\Cap_{p,\Omega}$-q.e. on }A} \\
     \label{eq:inf}
     &=& \inf\klg{\norm{\su}^p_{W^{1,p}(\Omega)}:\su\in\tW^{1,p}(\Omega),\tilde\su\geq 1\mbox{ $\Cap_{p,\Omega}$-q.e. on }A}.
  \end{eqnarray}
\end{theorem}

\begin{proof}
  Denote by $I$ the infimum on the right hand side of \eqref{eq:inf} and let $\su\in\cY_{p,\Omega}(A)$.
  Then by Theorem \ref{thm:ypo} $\tilde\su\geq 1$ $\Cap_{p,\Omega}$-q.e. on $A$. Hence
  $I\leq\norm{\su}^p_{W^{1,p}(\Omega)}$. Taking the infimum over all $\su\in\cY_{p,\Omega}(A)$ we
  get that $I\leq \Cap_{p,\Omega}(A)$. On the other hand, let $\su\in\tW^{1,p}(\Omega)$ be such that
  $\tilde\su\geq 1$ $\Cap_{p,\Omega}$-q.e. on $A$. Then by Lemma \ref{lem:approx} there exists
  $\su_n\in\cY_{p,\Omega}(A)$ such that $\su_n\to\su$ in $\tW^{1,p}(\Omega)$. Therefore
  $\Cap_{p,\Omega}(A)\leq\norm{\su_n}^p_{W^{1,p}(\Omega)}$ and passing to the limit as $n\to\infty$ gives
  $\Cap_{p,\Omega}(A)\leq\norm{\su}^p_{W^{1,p}(\Omega)}$. Now taking the infimum over all such
  $\su$ gives that $\Cap_{p,\Omega}(A)\leq I$ and hence we have equality.
\end{proof}

\begin{definition}\label{def:capextremal}
  A function $\su\in\tsW^{1,p}(\Omega)$ is called a/the
  {\em $\Cap_{p,\Omega}$-extremal for $A\subset\overline\Omega$} if $\su\geq 1$ $\Cap_{p,\Omega}$-q.e on $A$
  and $\norm{\su}_{W^{1,p}(\Omega)}^p=\Cap_{p,\Omega}(A)$.
\end{definition}

\begin{theorem}
  For every $A\subset\overline\Omega$ with $\Cap_{p,\Omega}(A)<\infty$ there exists a unique
  $\Cap_{p,\Omega}$-extremal $\se_A\in\tsW^{1,p}(\Omega)$.
  Moreover, $0\leq\se_A\leq 1$ $\Cap_{p,\Omega}$-q.e. on $\overline\Omega$ and
  $\se_A=1$ $\Cap_{p,\Omega}$-q.e. on $A$.
\end{theorem}

\begin{proof}
  Since $\Cap_{p,\Omega}(A)<\infty$ we have that $\overline\cY_{p,\Omega}(A)$ is a non-empty closed and
  convex subset of $\tW^{1,p}(\Omega)$. Let $(u_n)_n\subset\cY_{p,\Omega}(A)$ be such that
  $\norm{\su_n}^p_{W^{1,p}(\Omega)}\to\Cap_{p,\Omega}(A)$. Then the sequence $(\su_n)_n$ is
  bounded in the reflexive Banach space $\tW^{1,p}(\Omega)$ and hence, by possibly passing
  to a subsequence, weakly convergent to a function $\su\in\overline\cY_{p,\Omega}(A)$.
  Using the lower semi-continuity of the norm we get that
  \[ \norm{\su}^p_{W^{1,p}(\Omega)}\leq \liminf_n \norm{\su_n}^p_{W^{1,p}(\Omega)}=\Cap_{p,\Omega}(A).
  \]
  Since $\norm{\sv}^p_{W^{1,p}(\Omega)}\geq\Cap_{p,\Omega}(A)$ for all $\sv\in\cY_{p,\Omega}(A)$ we get
  that this equality remains true on $\overline\cY_{p,\Omega}(A)$ and hence
  $\norm{\su}^p_{W^{1,p}(\Omega)}=\Cap_{p,\Omega}(A)$. From Theorem \ref{thm:ypo} we get that
  $\se_A:=\tilde\su\geq 1$ $\Cap_{p,\Omega}$-q.e. on $A$. It remains to show the uniqueness.
  For this let $\sw\in\tsW^{1,p}(\Omega)$ be a $\Cap_{p,\Omega}$-extremal for $A$.
  Then $\sw\in\overline\cY_{p,\Omega}(A)$ by Theorem \ref{thm:ypo}. If $\sw\not=\se_A$ then by the uniform convexity
  of $\tW^{1,p}(\Omega)$ there exists $\sv\in\overline\cY_{p,\Omega}(A)$ with
  $\norm{\sv}^p_{W^{1,p}(\Omega)}<\Cap_{p,\Omega}(A)$ which is a contradiction.
  The additional claims follow as in Proposition \ref{prop:open-extr}.
\end{proof}

\begin{remark}\label{rem:inypo}
  The $\Cap_{p,\Omega}$-extremal for $A\subset\overline\Omega$ is the projection of
  $0$ onto $\overline\cY_{p,\Omega}(A)$.
\end{remark}
 
In the following we will use the convention that $\betrag{\xi}^{p-2}\xi:=0\in\IR^d$ if $\xi=0\in\IR^d$.

\begin{lemma}\label{lem:tool}
  Let $\su,\sv\in W^{1,p}(\Omega)$ and define $\sv_\varepsilon:=\su+\varepsilon\sv$ for $\varepsilon>0$. Then
  \[ \lim_{\varepsilon\to 0+} \varepsilon^{-1}\kle{\norm{\sv_\varepsilon}^p_{W^{1,p}(\Omega)}-\norm{\su}^p_{W^{1,p}(\Omega)}} =
     p\int_{\Omega} \betrag{\nabla\su}^{p-2}\nabla\su\nabla\sv+\betrag{\su}^{p-2}\su\sv.
  \]
\end{lemma}

\begin{proof}
  For $\xi,\theta\in\IR^d$ we have that (using the derivative with respect to $\varepsilon$)
  \[ \lim_{\varepsilon\to 0+} \varepsilon^{-1}\klr{\betrag{\xi+\varepsilon\theta}^p-\betrag{\xi}^p}=p\betrag{\xi}^{p-2}\xi\theta
  \]
  and thus by Lebesgue's Dominated Convergence Theorem
  \[ \lim_{\varepsilon\to 0+} \varepsilon^{-1}\kle{\norm{\sv_\varepsilon}^p_{W^{1,p}(\Omega)}-\norm{\su}^p_{W^{1,p}(\Omega)}} =
  \]
  \[ \lim_{\varepsilon\to 0+} \int_\Omega \varepsilon^{-1}\klr{\betrag{\su+\varepsilon\sv}^p-\betrag{\su}^p}+
             \varepsilon^{-1}\klr{\betrag{\nabla\su+\varepsilon\nabla\sv}^p-\betrag{\nabla\su}^p} =
  \]\[
     p\int_\Omega \betrag{\nabla\su}^{p-2}\nabla\su\nabla\sv+\betrag{\su}^{p-2}\su\sv.
  \]
\end{proof}

\begin{remark}
  Since the function $\varphi:\IR^d\to\IR$, $\xi\mapsto\betrag{\xi}^p$ is convex we get that
  \begin{equation}\label{eq:remark}
     \betrag{\xi'}^p-\betrag{\xi}^p \geq p\betrag{\xi}^{p-2}\xi(\xi'-\xi) \quad\mbox{ for all }\xi',\xi\in\IR^d.
  \end{equation}
\end{remark}

\begin{proposition}\label{prop:extchar}
  Let $A\subset\overline\Omega$. Then a function $\su\in\overline\cY_{p,\Omega}(A)$ is the
  $\Cap_{p,\Omega}$-extremal for $A$ if and only if
  \begin{equation}\label{eq:capextine}
     \int_\Omega \betrag{\nabla\su}^{p-2}\nabla\su\kle{\nabla\sv-\nabla\su}+\betrag{\su}^{p-2}\su\kle{\sv-\su} \geq 0\quad
     \mbox{ for all } \sv\in\overline\cY_{p,\Omega}(A).
  \end{equation}
\end{proposition}

\begin{proof}
  Let $\se_A$ be the $\Cap_{p,\Omega}$-extremal for $A$ and let $\sv\in\overline\cY_{p,\Omega}(A)$. For
  $\varepsilon>0$ we let $\sv_\varepsilon:=\se_A+\varepsilon(\sv-\se_A)$.
  Using that $\se_A=1$ $\Cap_{p,\Omega}$-q.e. on $A$ and Theorem \ref{thm:ypo} we get that
  $\sv_\varepsilon\in\overline\cY_{p,\Omega}(A)$ and hence
  \[ \varepsilon^{-1}\klr{\norm{\sv_\varepsilon}^p_{W^{1,p}(\Omega)}-\norm{\se_A}^p_{W^{1,p}(\Omega)}}\geq 0.
  \]
  Using Lemma \ref{lem:tool} we get \eqref{eq:capextine}. For the converse implication assume that
  $\su\in\overline\cY_{p,\Omega}(A)$ satisfies \eqref{eq:capextine}. Then by \eqref{eq:remark} we get that
  $\norm{\sv}^p_{W^{1,p}(\Omega)}-\norm{\su}^p_{W^{1,p}(\Omega)}\geq 0$ for all $\sv\in\overline\cY_{p,\Omega}(A)$
  and hence
  \[ \norm{\su}^p_{W^{1,p}(\Omega)}=\inf\klg{\norm{\sv}^p_{W^{1,p}(\Omega)}:\sv\in\cY_{p,\Omega}(A)}=\Cap_{p,\Omega}(A).
  \]
\end{proof}

\subsection{Potentials}

The aim of this subsection is to prove the existence and uniqueness of
$\tsW^{1,p}(\Omega)$-potentials for $\mu\in\tsW^{1,p}(\Omega)'$ and to
characterize them.

\begin{definition}
  Let $\mu\in\tsW^{1,p}(\Omega)'$ where $\tsW^{1,p}(\Omega)'$ is the topological dual of $\tsW^{1,p}(\Omega)$.
  Then $\su\in\tsW^{1,p}(\Omega)$ is called a/the {\em $\tsW^{1,p}(\Omega)$-potential of $\mu$} if
  $\su$ minimizes in $\tsW^{1,p}(\Omega)$ the mapping
  \begin{equation}\label{eq:mini}
     \tsW^{1,p}(\Omega)\to\IR,\quad \sv\mapsto \frac{1}{p}\norm{\sv}^p_{W^{1,p}(\Omega)}-\mu(\sv).
  \end{equation}
\end{definition}

\begin{theorem}\label{thm:uniexi}
  Let $\mu\in\tsW^{1,p}(\Omega)'$. Then there exists a unique $\su\in\tsW^{1,p}(\Omega)$ such that
  \begin{equation}\label{eq:minimizer}
     \frac{1}{p}\norm{\su}^p_{W^{1,p}(\Omega)}-\mu(\su) =
     \inf_{\sv\in\tsW^{1,p}(\Omega)} \frac{1}{p}\norm{\sv}^p_{W^{1,p}(\Omega)}-\mu(\sv).
  \end{equation}
  That is, for every $\mu\in\tsW^{1,p}(\Omega)'$ there exists a unique $\tsW^{1,p}(\Omega)$-potential of $\mu$.
\end{theorem}

\begin{proof}
  Let $(\su_n)_n\subset\tsW^{1,p}(\Omega)$ be a minimizing sequence of \eqref{eq:mini}.
  Since the infimum on the right hand side of \eqref{eq:minimizer} is less or equal to $0$
  (take $\sv=0\in\tsW^{1,p}(\Omega)$), we may assume that
  \[ \frac{1}{p}\norm{\su_n}^p_{W^{1,p}(\Omega)}\leq\mu(\su_n)\leq \norm{\mu}_{\tsW^{1,p}(\Omega)'}\norm{\su_n}\quad
     \mbox{ for all } n\in\IN.
  \]
  This shows that the sequence $(\su_n)_n$ is bounded in the reflexive Banach space $\tsW^{1,p}(\Omega)$
  and hence, by possibly passing to a subsequence, weakly convergent to a function $\su\in\tsW^{1,p}(\Omega)$.
  By using the lower semi-continuity of the norm in $\tsW^{1,p}(\Omega)$ we get that
  \[ \inf_{\sv\in\tsW^{1,p}(\Omega)} \frac{1}{p}\norm{\sv}^p_{\tsW^{1,p}(\Omega)}-\mu(\sv) =
     \lim_n \frac{1}{p}\norm{\su_n}^p_{\tsW^{1,p}(\Omega)}-\mu(\su_n) \geq
     \frac{1}{p}\norm{\su}^p_{\tsW^{1,p}(\Omega)}-\mu(\su).
  \]
  This shows that $\su$ is a minimizer and the existence is proved. The uniqueness follows from
  the strict convexity of $\norm{\cdot}^p_{W^{1,p}(\Omega)}$. In fact, assume that
  $\su_1,\su_2$ are two different minimizer of \eqref{eq:mini} and let $\su:=(\su_1+\su_2)/2$. Then 
  \begin{eqnarray*}
    \frac{1}{p}\norm{\su}^p_{\tsW^{1,p}(\Omega)}-\mu(\su) 
     &<& \frac{1}{2p}\kle{\norm{\su_1}^p_{\tsW^{1,p}(\Omega)}+\norm{\su_2}^p_{\tsW^{1,p}(\Omega)}}-\frac{\mu(\su_1+\su_2)}{2}\\
     &=& \inf_{\sv\in\tsW^{1,p}(\Omega)} \frac{1}{p}\norm{\sv}^p_{\tsW^{1,p}(\Omega)}-\mu(\sv),
  \end{eqnarray*}
  a contradiction.
\end{proof}

Theorem \ref{thm:uniexi} gives the existence and uniqueness for the $\tsW^{1,p}(\Omega)$-potential
for every $\mu\in\tsW^{1,p}(\Omega)'$. A characterization for this unique $\tsW^{1,p}(\Omega)$-potential
in terms of an integral equation is given by the following lemma.

\begin{lemma}\label{lem:wichtig}
  Let $\mu\in\tsW^{1,p}(\Omega)'$. Then $\su\in\tsW^{1,p}(\Omega)$ is the $\tsW^{1,p}(\Omega)$-potential
  of $\mu$ if and only if
  \begin{equation}\label{eq:equality}
     \int_\Omega \betrag{\su}^{p-2}\su\sv+\betrag{\nabla\su}^{p-2}\nabla\su\nabla\sv=\mu(\sv)\quad
     \mbox{ for all } \sv\in\tsW^{1,p}(\Omega).
  \end{equation}
\end{lemma}

\begin{proof}
  Let $\su\in\tsW^{1,p}(\Omega)$ be the $\tsW^{1,p}(\Omega)$-potential of $\mu$. For
  $\varepsilon>0$ and $\sv\in\tsW^{1,p}(\Omega)$ we let $\sv_\varepsilon:=\su+\varepsilon\sv\in\tsW^{1,p}(\Omega)$.
  Using the inequality
  \[ \frac{1}{p}\norm{\su}^p_{W^{1,p}(\Omega)} -\mu(\su)\leq
     \frac{1}{p}\norm{\sv_\varepsilon}^p_{W^{1,p}(\Omega)} -\mu(\sv_\varepsilon)
  \]
  we deduce that
  \[ \varepsilon^{-1}\frac{1}{p}\kle{\norm{\sv_\varepsilon}^p_{W^{1,p}(\Omega)}-\norm{\su}^p_{W^{1,p}(\Omega)}}\geq
     \mu\klr{\varepsilon^{-1}(\sv_\varepsilon-\su)}=\mu(\sv).
  \]
  For $\varepsilon\to 0+$ we get from Lemma \ref{lem:tool} that
  \[ \int_\Omega \betrag{\su}^{p-2}\su\sv+\betrag{\nabla\su}^{p-2}\nabla\su\nabla\sv\geq \mu(\sv).
  \]
  Replacing $\sv$ be $-\sv$ we get equality. Hence we proved that the $\tsW^{1,p}(\Omega)$-potential
  of $\mu$ satisfies \eqref{eq:equality}. To prove the sufficiency part let $\su\in\tsW^{1,p}(\Omega)$
  be such that $\su$ satisfies \eqref{eq:equality}. Then, by \eqref{eq:remark}, we get for
  $\sv\in\tsW^{1,p}(\Omega)$ that
  \[ \norm{\sv}^p_{W^{1,p}(\Omega)}-\norm{\su}^p_{W^{1,p}(\Omega)} \geq
     p\int_\Omega \betrag{\su}^{p-2}\su(\sv-\su)+\betrag{\nabla\su}^{p-2}\nabla\su(\nabla\sv-\nabla\su)=p\mu(\sv-\su),
  \]
  that is,
  \[ \frac{1}{p}\norm{\su}^p_{W^{1,p}(\Omega)}-\mu(\su) \leq \frac{1}{p}\norm{\sv}^p_{W^{1,p}(\Omega)}-\mu(\sv)
     \quad\mbox{ for all }\sv\in\tsW^{1,p}(\Omega).
  \]
\end{proof}

\begin{theorem}\label{thm:capmeasure}
  Let $\mu\in\tsW^{1,p}(\Omega)'$ and $\su\in\tsW^{1,p}(\Omega)$ be the $\tsW^{1,p}(\Omega)$-potential of $\mu$.
  Then
  \[ \norm{\mu}^{p'}_{\tsW^{1,p}(\Omega)'}=\mu(\su)=\norm{\su}^p_{W^{1,p}(\Omega)}\qquad
     \mbox{ with } 1/p+1/p'=1.
  \]
\end{theorem}

\begin{proof}
  It follows immediately from Lemma \ref{lem:wichtig} that $\mu(\su)=\norm{\su}^p_{W^{1,p}(\Omega)}$.
  To prove the remaining equality let $\sv\in\tsW^{1,p}(\Omega)$ be such that
  $\norm{\sv}_{W^{1,p}(\Omega)}\leq\norm{\su}_{W^{1,p}(\Omega)}$. Using the definition of the
  $\tsW^{1,p}(\Omega)$-potential we get that
  \begin{eqnarray*}
     \mu(\sv)
        &=& \kle{\mu(\sv)-\frac{1}{p}\norm{\sv}^p_{W^{1,p}(\Omega)}}+\frac{1}{p}\norm{\sv}^p_{W^{1,p}(\Omega)} \\
        &\leq& \kle{\mu(\su)-\frac{1}{p}\norm{\su}^p_{W^{1,p}(\Omega)}}+\frac{1}{p}\norm{\su}^p_{W^{1,p}(\Omega)} = \mu(\su).
  \end{eqnarray*}
  It follows from Lemma \ref{lem:wichtig} that if $\su=0$ then $\mu=0$ and the assertion is trivial.
  So we may assume that $\su\not=0$. Using the linearity of $\mu$ we get that for all $\sw\in\tsW^{1,p}(\Omega)$
  with $\norm{\sw}_{W^{1,p}(\Omega)}\leq 1$
  \[ \mu(\sw) \leq \mu(u/\norm{u}_{W^{1,p}(\Omega)})
  \]
  and hence
  \[ \norm{\mu}^{p'}_{\tsW^{1,p}(\Omega)'} = \mu(\su/\norm{\su}_{W^{1,p}(\Omega)})^{p'} 
     = \mu(\su)^{p'}\norm{\su}^{-p'}_{W^{1,p}(\Omega)} = \norm{\su}^{(p-1)p'}_{W^{1,p}(\Omega)}=\norm{\su}^p_{W^{1,p}(\Omega)}.
  \]
\end{proof}

\subsection{Capacitary Measures}

The aim of this subsection is to prove the existence and uniqueness for
$\Cap_{p,\Omega}$-measures and to characterize positive functionals
in $\tsW^{1,p}(\Omega)'$.

\begin{definition}
  A Borel measure $\mu$ on a Hausdorff space $\topo$ is called a {\em Radon measure}
  if $\mu$ is inner regular and locally finite.
\end{definition}

\begin{definition}\label{def:radonin}
  We say that a Radon measure $\mu$ on $\overline\Omega$ belongs to $\tsW^{1,p}(\Omega)'$ if
  \[ \sup\klg{\int_{\overline\Omega} v\;d\mu:v\in\tsW^{1,p}(\Omega)\cap C_c(\overline\Omega),
              \norm{v}_{W^{1,p}(\Omega)}\leq 1}<\infty.
  \]
  In this case we can extend the functional $v\mapsto \int_{\overline\Omega} v\;d\mu$ from
  $\tsW^{1,p}(\Omega)\cap C_c(\overline\Omega)$ in a unique way to a continuous functional
  $\tilde\mu$ on $\tsW^{1,p}(\Omega)$, i.e. $\tilde\mu\in\tsW^{1,p}(\Omega)'$.
  For $\sv\in\tsW^{1,p}(\Omega)$ we simply write $\mu(\sv)$ instead of $\tilde\mu(\sv)$ but
  still have this definition of $\mu(\sv)$ in mind.
\end{definition}

\begin{definition}
  Let $A\subset\overline\Omega$ be such that $\Cap_{p,\Omega}(A)<\infty$ and let
  $\mu$ be a Radon measure in $\tsW^{1,p}(\Omega)'$. Then $\mu$ is said to be a
  {\em $\Cap_{p,\Omega}$-measure for $A$} if the $\tsW^{1,p}(\Omega)$-potential of $\mu$
  is the $\Cap_{p,\Omega}$-extremal for $A$.
\end{definition}

\begin{definition}\label{def:energy}
  Let $\mu\in\tsW^{1,p}(\Omega)'$. Then the {\em $\Cap_{p,\Omega}$-energy} of $\mu$ on $\overline\Omega$
  is defined by
  \[ \sE(\mu):=\sE_p(\mu,\Omega):=\norm{\mu}^{p'}_{\tsW^{1,p}(\Omega)'}.
  \]
\end{definition}

\begin{lemma}\label{lem:ext}
  Let $\psi\in\tsW^{1,p}(\Omega)'$ be such that $\psi(w)\geq 0$ for all non-negative
  $w\in S:=\tsW^{1,p}(\Omega)\cap C_c(\overline\Omega)$. Then $\psi|_S$
  can be extended to a positive functional $\tilde\psi$ on $C_c(\overline\Omega)$.
\end{lemma}

\begin{proof}
  Let $u\in C_c(\overline\Omega)$ and $K:=\supp(u)\subset\overline\Omega$.
  Fix a non-negative test function $\varphi\in\cD(\IR^N)$ such that $\varphi\equiv 1$ on $K$
  and $0\leq\varphi\leq 1$ on $\IR^N$. Let $u_n\in\cD(\IR^N)$
  be such that $u_n$ converges to $u$ uniformly on $\overline\Omega$.
  Since $\betrag{\psi(v\varphi)}\leq \psi(\varphi)\norm{v}_{\infty}$
  for all $v\in S$ we get that $\tilde\psi(u):=\lim_n\psi(\varphi u_n)$
  exists. Note that the definition of $\tilde\psi$ does not depend on the
  sequence $(u_n)_n$ and $\varphi$ and that $\tilde\psi$ is linear and positive.
\end{proof}

\begin{lemma}\label{lem:radon}
  Let $\psi\in\tsW^{1,p}(\Omega)'$ be such that $\psi(w)\geq 0$ for all non-negative
  $w\in S:=\tsW^{1,p}(\Omega)\cap C_c(\overline\Omega)$. Then there exists a unique Radon measure
  $\mu$ on $\overline\Omega$ such that
  \begin{equation}\label{eq:radon}
     \psi(w) = \int_{\overline\Omega} w\;d\mu\quad\mbox{ for all } w\in S=\tsW^{1,p}(\Omega)\cap C_c(\overline\Omega).
  \end{equation}
\end{lemma}

\begin{proof}
  Let $\tilde\psi$ be given from Lemma \ref{lem:ext}. Then we get from the Riesz-Markov Theorem
  (see Royden \cite[Theorem 13.23]{royden:88:ran}) that there exists a unique Radon measure $\mu$ on $\overline\Omega$
  such that $\tilde\psi(w)=\int w\;d\mu$ for all $w\in C_c(\overline\Omega)$. In particular,
  $\psi(w)=\tilde\psi(w)=\int w\;d\mu$ for all $w\in S$. Let $\mu'$ be a Radon measure on $\overline\Omega$ satisfying
  \eqref{eq:radon}. Using the density of $S$ in $C_c(\overline\Omega)$ we get that
  \[ \int_{\overline\Omega} w\;d\mu = \int_{\overline\Omega} w\;d\mu'\;\quad\mbox{ for all } w\in C_c(\overline\Omega)
  \]
  and hence by the Riesz-Markov Representation Theorem that $\mu=\mu'$.
\end{proof}

\begin{theorem}
  Let $A\subset\overline\Omega$ be such that $\Cap_{p,\Omega}(A)<\infty$.
  Then there exists a unique $\Cap_{p,\Omega}$-measure $\mu\in\tsW^{1,p}(\Omega)'$ for $A$.
  Moreover, if $\su\in\tsW^{1,p}(\Omega)$ is the $\Cap_{p,\Omega}$-extremal for $A$ then
  \[ \Cap_{p,\Omega}(A)=\norm{\su}^p_{W^{1,p}(\Omega)}=\mu(\su)=\sE_p(\mu,\Omega).
  \]
\end{theorem}

\begin{proof}
  By Proposition \ref{prop:extchar} and Remark \ref{rem:inypo} the $\Cap_{p,\Omega}$-extremal $\su$ for $A$ is
  in $\overline\cY_{p,\Omega}(A)$ and satisfies
  \begin{equation}\label{eq:potchar}
     \int_{\Omega} \betrag{\nabla\su}^{p-2}\nabla\su\kle{\nabla\sv-\nabla\su}+\betrag{\su}^{p-2}\su\kle{\sv-\su} \geq 0
     \quad\mbox{ for all } \sv\in\overline\cY_{p,\Omega}(A).
  \end{equation}
  We define $\psi\in\tsW^{1,p}(\Omega)'$ by
  \[ \psi(\sw):=\int_{\Omega} \betrag{\nabla\su}^{p-2}\nabla\su\nabla\sw+\betrag{\su}^{p-2}\su\sw.
  \]
  If $w\in S:=\tsW^{1,p}(\Omega)\cap C_c(\overline\Omega)$ is non-negative, then
  $\sv:=\su+w\in\overline\cY_{p,\Omega}(A)$ and hence we get from \eqref{eq:potchar}
  \begin{eqnarray*}
     \psi(w) 
       &=& \int_{\Omega} \betrag{\nabla\su}^{p-2}\nabla\su\nabla w+\betrag{\su}^{p-2}\su w\\
       &=& \int_{\Omega} \betrag{\nabla\su}^{p-2}\nabla\su\kle{\nabla\sv-\nabla\su}+\betrag{\su}^{p-2}\su\kle{\sv-\su}\geq 0.
  \end{eqnarray*}
  By Lemma \ref{lem:radon} we get that there exists a unique Radon measure $\mu$ on $\overline\Omega$ such that
  \[ \psi(w) = \int_{\overline\Omega} w\;d\mu
     \quad\mbox{ for all } w\in\tsW^{1,p}(\Omega)\cap C_c(\overline\Omega).
  \]
  It follows from the continuity of $\psi$ and Definition \ref{def:radonin} that $\psi(\sw)=\mu(\sw)$ for all
  $\sw\in\tsW^{1,p}(\Omega)$. From Lemma \ref{lem:wichtig} we get that $\su$ is the $\tsW^{1,p}(\Omega)$-potential
  of $\mu$ and hence $\mu$ is a $\Cap_{p,\Omega}$-measure for $A$ by definition. If $\mu'$ is a Radon measure
  in $\tsW^{1,p}(\Omega)'$ such that $\su$ is the $\Cap_{p,\Omega}$-potential of $\mu'$, then by
  Lemma \ref{lem:wichtig} $\mu'(w)=\mu(w)$ for all $w\in S\subset\tsW^{1,p}(\Omega)$,
  from which we deduce the uniqueness of the $\Cap_{p,\Omega}$-measure for $A$.
  The stated equality follows from Theorem \ref{thm:capmeasure} and Definitions \ref{def:energy}
  and \ref{def:capextremal}.
\end{proof}

\begin{theorem}\label{thm:abscont}
  Let $\mu\in\tsW^{1,p}(\Omega)'$ be a Radon measure and let $A\subset\overline\Omega$ be $\mu$-measurable. Then
  \begin{equation}\label{eq:claim}
     \mu(A)^p \leq \sE(\mu)^{p-1}\Cap_{p,\Omega}(A).
  \end{equation}
\end{theorem}

\begin{proof}
  First assume that $A$ is compact. Then for any non-negative $v\in\tsW^{1,p}(\Omega)\cap C_c(\overline\Omega)$
  with $v\geq 1$ on $A$ we have that
  \[ \mu(A) \leq \int_{\overline\Omega} v\;d\mu = \mu(v)\leq\norm{\mu}_{\tsW^{1,p}(\Omega)'}\norm{v}_{W^{1,p}(\Omega)}.
  \]
  Taking the infimum over all such $v$ we get by Proposition \ref{prop:cont} that
  \[ \mu(A) \leq \sE_p(\mu,\Omega)^{1/p'}\Cap_{p,\Omega}(A)^{1/p}
  \]
  which is equivalent to \eqref{eq:claim}. If $A$ is an open set in $\overline\Omega$ then we
  consider an increasing sequence of compact sets $(A_n)_n$ such that $A=\bigcup_n A_n$. Then
  \[ \mu(A) = \lim_n \mu(A_n) \leq \lim_n \sE_p(\mu,\Omega)^{1/p'}\Cap_{p,\Omega}(A_n)^{1/p} =
          \sE_p(\mu,\Omega)^{1/p'}\Cap_{p,\Omega}(A)^{1/p}.
  \]
  Finally let $A$ be an arbitrary $\mu$-measurable set, then (since every Radon measure on $\overline\Omega$
  is automatically outer regular -- see Royden \cite[Proposition 13.14]{royden:88:ran}) we get
  \begin{eqnarray*}
    \mu(A) &=& \inf\klg{\mu(O):O\supset A\mbox{ open in }\overline\Omega} \\
        &\leq& \inf\klg{\sE_p(\mu,\Omega)^{1/p'}\Cap_{p,\Omega}(O)^{1/p}:O\supset A\mbox{ open in }\overline\Omega} \\
        &=& \sE_p(\mu,\Omega)^{1/p'}\Cap_{p,\Omega}(A).
  \end{eqnarray*}
\end{proof}

\begin{corollary}
  If $\mu\in\tsW^{1,p}(\Omega)'$ is a Radon measure, then $\mu(A)=0$ for every
  $\Cap_{p,\Omega}$-polar set. That is, $\mu$ is absolutely continuous with respect to $\Cap_{p,\Omega}$.
\end{corollary}

\begin{theorem}
  Let $\mu\in\tsW^{1,p}(\Omega)'$ be a Radon measure.
  Then $\tsW^{1,p}(\Omega)\subset L^1(\overline\Omega,\mu)$ and
  \[ \mu(\su) = \int_{\overline\Omega} \su\;d\mu\qquad\text{ for all }\su\in\tsW^{1,p}(\Omega).
  \]
\end{theorem}

\begin{remark}
  Note that $\tsW^{1,p}(\Omega)$ consists only of $\Cap_{p,\Omega}$-quasi continuous functions and so there is
  no need to pass from a function $\sv\in\tW^{1,p}(\Omega)$ to its $\Cap_{p,\Omega}$-quasi continuous
  representative $\tilde\sv:\overline\Omega\to\IR$.
\end{remark}

\begin{proof}
  Let $\su\in\tsW^{1,p}(\Omega)$ and let $(u_k)_k$ be a sequence in $\tsW^{1,p}(\Omega)\cap C_c(\overline\Omega)$
  such that $(u_k)_k$ converges to $\su$ in $\tsW^{1,p}(\Omega)$. By possibly passing to a subsequence we may
  assume that $\norm{u_k-\su}_{W^{1,p}(\Omega)}\leq 4^{-k}$ for all $k\in\IN$.
  Let $w:=\sum_{k\in\IN} \betrag{u_{k+1}-u_k}$. Then
  \begin{eqnarray*}
     \int_{\overline\Omega} w\;d\mu 
       &=& \sum_{k=1}^\infty \int_{\overline\Omega} \betrag{u_{k+1}-u_k}\;d\mu =
           \sum_{k=1}^\infty \norm{u_{k+1}-u_k}_{L^1(\overline\Omega,\mu)} \\
       &=& \sum_{k=1}^\infty \mu\klr{\betrag{u_{k+1}-u_k}} \leq
           \sum_{k=1}^\infty \norm{\mu}_{\tsW^{1,p}(\Omega)'}\norm{u_{k+1}-u_k}_{W^{1,p}(\Omega)}<\infty.
  \end{eqnarray*}
  First we show that $(u_k(x))_k$ is a Cauchy sequence for $\mu$-a.e.
  $x\in\overline\Omega$. Let $x\in\overline\Omega$ be fixed. If $(u_k(x))_k$ is not a Cauchy sequence
  then $w(x)=\infty$ and since $w\in L^1(\overline\Omega,\mu)$ we get that the set where
  $w=\infty$ is a $\mu$-nullset, that is, $(u_k(x))_k$ is a Cauchy sequence in $\IR$ for $\mu$-almost
  every $x\in\overline\Omega$. It also follows immediately from the equation above that
  $(u_k)_k$ is a Cauchy sequence in $L^1(\overline\Omega,\mu)$ and hence convergent to a function
  $\sv\in L^1(\overline\Omega,\mu)$. Moreover, $u_k$ converges $\mu$-a.e. on $\overline\Omega$ to $\sv$.
  On the other hand, it follows from Theorem \ref{thm:convsubse} that, by possibly passing to a subsequence,
  $u_k$ converges $\Cap_{p,\Omega}$-quasi everywhere (and hence $\mu$-a.e.) on $\overline\Omega$ to $\su$.
  Therefore $\su=\sv$ $\mu$-a.e. on $\overline\Omega$ and
  \[ \int_{\overline\Omega} \su\;d\mu = \int_{\overline\Omega} \sv\;d\mu =
     \lim_k \int_{\overline\Omega} u_k\;d\mu = \lim_k \mu(u_k) = \mu(\su).
  \]
  In particular $\su\in L^1(\overline\Omega,\mu)$.
\end{proof}

\section{An application to Sobolev spaces}

In this section we give an application of the relative capacity, namely to decide
if a given function $\su$ lies in $W^{1,p}_0(\Omega)$ or not.
Here $W^{1,p}_0(\Omega)$ is the closure of $\cD(\Omega)$ in $W^{1,p}(\Omega)$ and hence a closed
subspace of $\tilde W^{1,p}(\Omega)$. For an open set $\Omega\subset\IR^N$ the following inclusion holds.
\[ W^{1,p}(\Omega)\cap C_0(\Omega)\subset W^{1,p}_0(\Omega),\quad 1\leq p<\infty.
\]
Here $C_0(\Omega)$ is defined to be the space of all continuous functions $u:\Omega\to\IR$ such
that for all $\varepsilon>0$ there exists a compact set $K\subset\Omega$ such that
$\betrag{u}\leq\varepsilon$ for all $x\in\Omega\setminus K$. To prove this inclusion let
$u\in W^{1,p}(\Omega)\cap C_0(\Omega)$. Since $u=u^+-u^-$ we may assume without loss of
generality that $u\geq 0$. For $k\in\IN$ let $u_k:=(u-1/k)^+$. Then $u_k$ has compact support
in $\Omega$. Using a mollification argument
($\cD(\Omega)\ni u_{k,n}:=\rho_n\star u_k\to u_k$ in $W^{1,p}(\Omega)$) we see that $u_k\in W^{1,p}_0(\Omega)$
and since $u_k\to u$ in $W^{1,p}(\Omega)$ the claim is proved.

\begin{lemma}\label{lem:w1p0}
  Assume that $\Omega\subset\IR^N$ is a bounded and non-empty open set and $\su\in W^{1,p}(\Omega)$.
  If there exists $u\in\su$ such that for all $z\in\partial\Omega$ the limit
  $\lim_{\Omega\ni x\to z} u(x)$ exists and is equal to $0$, then $\su\in W^{1,p}_0(\Omega)$.
\end{lemma}

\begin{proof}
  By splitting $u=u^+-u^-$ we may assume that $u$ is a non-negative function.
  Now let $m\in\IN$ be fixed. By assumption, for $z\in\partial\Omega$, there exists
  $\delta_z>0$ such that $0\leq u(x)\leq 1/m$ for all $x\in B_{\IR^N}(z,\delta_z)\cap\Omega$.
  Since the boundary $\partial\Omega$ of $\Omega$ is compact, there exist
  $z_1,\dots,z_{n_0}\in\partial\Omega$ such that
  \[ \partial\Omega\subset\bigcup_{k=1}^{n_0} B_{\IR^N}(z_k,\delta_{z_k})=:O\subset\IR^N.
  \]
  Then $u_m:=(u-1/m)^+\in W^{1,p}(\Omega)$ and $u_m=0$ outside the compact
  set $K:=\overline\Omega\setminus O\subset\Omega$. Hence by the mollification argument
  described at the beginning of this section we get that $u_m\in W^{1,p}_0(\Omega)$.
  Since $u_m\to u$ in $W^{1,p}(\Omega)$ we deduce that $u\in W^{1,p}_0(\Omega)$.
\end{proof}

\begin{theorem}
  Let $\Omega\subset\IR^N$ be an open and non-empty set and $p\in (1,\infty)$. Then
  \begin{equation}\label{eq:w1p0}
     W^{1,p}_0(\Omega) = \klg{ \su\in\tW^{1,p}(\Omega):
                               \tilde\su=0\;  \mbox{ $\Cap_{p,\Omega}$-q.e. on $\partial\Omega$} }. 
  \end{equation}
\end{theorem}

\begin{proof}
  Let $D_0^{1,p}(\Omega)$ denote the right hand side of \eqref{eq:w1p0}.
  First we show that $W^{1,p}_0(\Omega)\subset D_0^{1,p}(\Omega)$.
  Let $\su\in W^{1,p}_0(\Omega)\subset\tW^{1,p}(\Omega)$. Then there exists
  a sequence of test functions $u_n\in\cD(\Omega)$ such that $u_n\to \su$ in $\tW^{1,p}(\Omega)$.
  By possibly passing to a subsequence (see Theorem \ref{thm:convsubse})
  we get that $(u_n)_n$ converges $\Cap_{p,\Omega}$-quasi everywhere to $\tilde\su$ and hence
  $\tilde\su=0$ $\Cap_{p,\Omega}$-quasi everywhere on $\partial\Omega$, that is, $\su\in D_0^{1,p}(\Omega)$.

  We show that $D_0^{1,p}(\Omega)\subset W^{1,p}_0(\Omega)$.
  Assume for the moment that $\Omega$ is bounded and let $\su\in D^{1,p}_0(\Omega)\cap L^\infty(\Omega)$
  be non-negative. Then there exists a sequence $(u_n)_n$ in
  $\tW^{1,p}(\Omega)\cap C_c(\overline\Omega)$ which converges to $\su$ in $\tW^{1,p}(\Omega)$.
  Since $(u_n\vee 0)\wedge\norm{\su}_\infty$ converges also to $\su$ in $\tW^{1,p}(\Omega)$ we may
  assume that $0\leq u_n\leq \norm{\su}_\infty$. Let $u\in\tilde\su$ be fixed. 
  By possibly passing to a subsequence (see Theorem \ref{thm:convsubse})
  we have that for each $m\in\IN$ there exists an open set $G_m$ in $\overline\Omega$ such that
  $\Cap_{p,\Omega}(G_m)\leq 1/m$ and $u_n\to u$ uniformly on $\overline\Omega\setminus G_m$.
  Hence there exists $n_0=n_0(m)$ such that $\betrag{u_{n_0}-u}\leq 1/m$ everywhere on $\overline\Omega\setminus G_m$
  and $\norm{\su_{n_0}-\su}_{W^{1,p}(\Omega)}\leq 1/m$. Let $U_m$ be an open set in $\overline\Omega$ such that
  $\Cap_{p,\Omega}(U_m)\leq 1/m$ and $u=0$ everywhere on $\partial\Omega\setminus U_m$.
  Consequently, $\betrag{u_{n_0}}\leq 1/m$ everywhere on $\partial\Omega\setminus O_m$ where $O_m:=G_m\cup U_m$.
  Let $\se_m\in\tsW^{1,p}(\Omega)$ be the $\Cap_{p,\Omega}$-extremal for $O_m$ and fix $e_m\in\se_m$.
  By changing $e_m$ on a $\Cap_{p,\Omega}$-polar set we may assume that $e_m\equiv 1$ everywhere on $O_m$
  and $0\leq e_m\leq 1$ everywhere on $\overline\Omega$.
  Let $w_m:=(u_{n_0}-1/m)^+$. Then we have that $v_m:=w_m(1-e_m)\in W^{1,p}_0(\Omega)$. In fact,
  for $z\in\partial\Omega$ we take a sequence $(x_n)_n$ in $\Omega$ which converges to $z$.
  If $z\in\partial\Omega\setminus O_m$ then
  \[ 0\leq v_m(x_n)\leq w_m(x_n)\to w_m(z)=(u_{n_0}(z)-1/m)^+=0.
  \]
  If $z\in\partial\Omega\cap O_m$ then there exists $k_0$ such that $x_n\in O_m$ for all $n\geq k_0$. Hence
  \[ 0\leq v_m(x_n)=w_m(x_n)(1-e_m(x_n))=0\quad\mbox{ for all } n\geq k_0.
  \]
  This shows that $\lim_{\Omega\ni x\to z} v_m(x)$ exists and is $0$ for all $z\in\partial\Omega$.
  It follows from Lemma \ref{lem:w1p0} that $v_m\in W^{1,p}_0(\Omega)$.
  Next we show that $v_m\to\su$ in $L^p(\Omega)$.
  \begin{eqnarray*}
    \norm{\su-v_m}_{L^p(\Omega)}
      & \leq & \norm{\su-u_{n_0}}_{L^p(\Omega)} +
               \norm{u_{n_0}-w_m}_{L^p(\Omega)} +
               \norm{w_m-v_m}_{L^p(\Omega)} \\
      & \leq & 1/m + \betrag{\Omega}^{1/p}/m +
               (2/m)^{1/p} \norm{w_m}_{L^\infty(\Omega)} \\
      & \leq & 1/m \kle{1+\betrag{\Omega}^{1/p}} + (2/m)^{1/p}\norm{u}_{L^\infty(\Omega)}.
  \end{eqnarray*}
  Next we estimate $\norm{D_j v_m}_{L^p(\Omega)}$ for $j\in\klg{1,\dots,N}$.
  \begin{eqnarray*}
    \norm{D_j v_m}_{L^p(\Omega)}
      & \leq & \norm{D_j w_m (1-e_m)}_{L^p(\Omega)} +
               \norm{w_m D_j e_m}_{L^p(\Omega)} \\
      & \leq & \norm{w_m}_{W^{1,p}(\Omega)} +
               \norm{w_m}_{L^\infty(\Omega)}\norm{e_m}_{W^{1,p}(\Omega)} \\
      & \leq & \norm{u_{n_0}}_{W^{1,p}(\Omega)} +
               \norm{u_{n_0}}_{L^\infty(\Omega)}\cdot \norm{e_m}_{W^{1,p}(\Omega)} \\
      & \leq & \norm{\su}_{W^{1,p}(\Omega)} + 1/m + \norm{\su}_{L^\infty(\Omega)}\cdot (2/m)^{1/p}
  \end{eqnarray*}
  This shows that the sequence $(v_m)_m$ is bounded in the reflexive Banach space
  $W^{1,p}_0(\Omega)$ and hence there exists a weakly convergent subsequence
  $(v_{m_k})_k$. Using that $W^{1,p}_0(\Omega)$ is closed for the weak topology,
  we get that the weak limit $\wlim_k v_{m_k}\in W^{1,p}_0(\Omega)$. On the other hand,
  since $(v_{m_k})_k$ converges to $\su$ in $L^p(\Omega)$ as $k\to\infty$,
  it follows that $\su=\wlim_k v_{m_k}\in W^{1,p}_0(\Omega)$. If $\su\in D^{1,p}_0(\Omega)\cap L^\infty(\Omega)$
  is not non-negative, we get by what we proved already that $\su^+$ and $\su^-$ are in $W^{1,p}_0(\Omega)$
  and hence $\su\in W^{1,p}_0(\Omega)$. If $\su\in D^{1,p}_0(\Omega)$ let
  $\sg_k:=(\su\wedge k)\vee(-k)\in D_0^{1,p}(\Omega)\cap L^\infty(\Omega)$. Since $\sg_k\to\su$
  in $W^{1,p}(\Omega)$ and $\sg_k\in W^{1,p}_0(\Omega)$ by the arguments above we get that $\su\in W^{1,p}_0(\Omega)$.
  If $\Omega$ is unbounded then we choose $\psi\in\cD(B(0,1))$ such that
  $\psi\equiv 1$ on a neighbourhood of $\klg{0}$. Let $\sv_n(x):=\su(x)\cdot\psi(x/n)$
  where $\su\in D^{1,p}_0(\Omega)$. Then by the above arguments we have
  \[ \sv_n\in D^{1,p}_0(\Omega\cap \cB(0,n))\subset W^{1,p}_0(\Omega\cap B(0,n)) \subset W^{1,p}_0(\Omega).
  \]
  Using that $\sv_n\to\su$ in $W^{1,p}(\Omega)$ the proof is finished.
\end{proof}

To finish this section and the article we mention two further characterizations of $W^{1,p}_0(\Omega)$.
The original proof of the following result is due to Havin \cite{havin:68:aim} and
Bagby \cite{bagby:72:qut}, an alternative proof is given by Hedberg \cite{hedberg:72:nlp}.

\begin{theorem}
  Let $1<p<\infty$, $\Omega\subset\IR^N$ an open set and let $\su\in W^{1,p}(\IR^N)$.
  Then $\su\in W^{1,p}_0(\Omega)$ if and only if
  \[ \lim_{r\to 0} r^{-N}\int_{B(x,r)} \betrag{u(y)}\;dy = 0
  \]
  for $Cap_{p}$-q.e. $x\in\IR^N\setminus\Omega$. 
\end{theorem}

The following characterization was recently proved by David Swanson and
William P. Ziemer. The main difference to the previous theorem is that the function
$\su$ is not assumed to belong to the space $W^{1,p}(\IR^N)$.

\begin{theorem}{D.~Swanson and W.~P.~Ziemer\bf\cite[Theorem 2.2]{swanson:99:sfi}}\label{thm:ziemer}.
  Let $\su\in W^{1,p}(\Omega)$. If
  \[ \lim_{r\to 0} r^{-N} \int_{\cB(x,r)\cap\Omega} \betrag{\su(y)}\;dy = 0
  \]
  for $\Cap_p$-quasi every $x\in\partial\Omega$, then $\su\in W^{1,p}_0(\Omega)$.
\end{theorem}

\bibliography{biblio2}

\providecommand{\mathbb}[1]{\mathbf{#1}}\providecommand{\cprime}{$'$}
\providecommand{\bysame}{\leavevmode\hbox to3em{\hrulefill}\thinspace}
\providecommand{\MR}{\relax\ifhmode\unskip\space\fi MR }
\providecommand{\MRhref}[2]{%
  \href{http://www.ams.org/mathscinet-getitem?mr=#1}{#2}
}
\providecommand{\href}[2]{#2}
\begin{thebibliography}{10}

\bibitem{adams:96:fsp}
David~R. Adams and Lars~Inge Hedberg, \emph{Function spaces and potential
  theory}, Grundlehren der mathematischen Wissenschaften, vol. 314,
  Springer-Verlag, Berlin, 1996. \MR{97j:46024}

\bibitem{alt:99:lfa}
Hans~Wilhelm Alt, \emph{{Linear functional analysis. An application oriented
  introduction. (Lineare Funktionalanalysis. Eine anwendungsorientierte
  Einf\"uhrung.) 3., vollst. \"uberarb. und erw. Auflage.}},
  {Springer-Lehrbuch. Berlin: Springer. xiii, 415 S.}, 1999 (German).

\bibitem{arendt:03:lrb}
Wolfgang Arendt and Mahamadi Warma, \emph{The {L}aplacian with {R}obin boundary
  conditions on arbitrary domains}, Potential Anal. \textbf{19} (2003), no.~4,
  341--363. \MR{1 988 110}

\bibitem{bagby:72:qut}
Thomas Bagby, \emph{{Quasi topologies and rational approximation.}}, J. Funct.
  Anal. \textbf{10} (1972), 259--268.

\bibitem{biegert:05:epvd}
Markus Biegert, \emph{Elliptic problems on varying domains}, Dissertation,
  Logos Verlag, Berlin, 2005.

\bibitem{bouleau:91:df}
Nicolas Bouleau and Francis Hirsch, \emph{Dirichlet forms and analysis on
  {W}iener space}, de Gruyter Studies in Mathematics, vol.~14, Walter de
  Gruyter \& Co., Berlin, 1991. \MR{MR1133391 (93e:60107)}

\bibitem{choquet:54:toc}
Gustave Choquet, \emph{{Theory of capacities.}}, Ann. Inst. Fourier (Grenoble)
  \textbf{5} (1954), 131--295.

\bibitem{dautray:88:man2}
Robert Dautray and Jacques-Louis Lions, \emph{Mathematical analysis and
  numerical methods for science and technology. {V}ol. 2}, Springer-Verlag,
  Berlin, 1988. \MR{89m:00001}

\bibitem{doob:01:cpt}
Joseph~L. Doob, \emph{{Classical potential theory and its probabilistic
  counterpart. Reprint of the 1984 edition.}}, {Classics in Mathematics.
  Berlin: Springer. xxiii, 846 p.}, 2001.

\bibitem{edmunds:87:std}
D.~E. Edmunds and W.~D. Evans, \emph{Spectral theory and differential
  operators}, Clarendon Press, Oxford, 1987. \MR{89b:47001}

\bibitem{evans:92:mtf}
Lawrence~C. Evans and Ronald~F. Gariepy, \emph{Measure theory and fine
  properties of functions}, Studies in Advanced Mathematics, CRC Press, Boca
  Raton, FL, 1992. \MR{MR1158660 (93f:28001)}

\bibitem{gilbarg:01:epd}
David Gilbarg and Neil~S. Trudinger, \emph{Elliptic partial differential
  equations of second order}, Classics in Mathematics, Springer-Verlag, Berlin,
  2001, Reprint of the 1998 edition. \MR{2001k:35004}

\bibitem{koskela:08:see}
Piotr Haj{\l}asz, Pekka Koskela, and Heli Tuominen, \emph{Sobolev embeddings,
  extensions and measure density condition}, J. Funct. Anal. \textbf{254}
  (2008), no.~5, 1217--1234. \MR{MR2386936}

\bibitem{havin:68:aim}
Victor~P. Havin, \emph{Approximation in the mean by analytic functions}, Dokl.
  Akad. Nauk SSSR 9 (1968), 245--248.

\bibitem{hedberg:72:nlp}
Lars~Inge Hedberg, \emph{{Non-linear potentials and approximation in the mean
  by analytic functions.}}, Math. Z. \textbf{129} (1972), 299--319.

\bibitem{heinonen:93:npt}
Juha Heinonen, Tero Kilpel{\"a}inen, and Olli Martio, \emph{Nonlinear potential
  theory of degenerate elliptic equations}, Oxford Mathematical Monographs,
  Clarendon Press, New York, 1993. \MR{94e:31003}

\bibitem{heuser:92:fa}
Harro Heuser, \emph{{Functional Analysis. Theory and Applications.
  (Funktionalanalysis. Theorie und Anwendung. Mit 766 Aufgaben, zum Teil mit
  L\"osungen und zahlreichen Beispielen.) 3., durchges. Aufl.}}, {Mathematische
  Leitf\"aden. Stuttgart: Teubner. 696 S. }, 1992 (German).

\bibitem{ziemer:97:fr}
Jan Mal{\'y} and William~P. Ziemer, \emph{Fine regularity of solutions of
  elliptic partial differential equations}, Mathematical Surveys and
  Monographs, vol.~51, American Mathematical Society, Providence, RI, 1997.
  \MR{MR1461542 (98h:35080)}

\bibitem{mazya:85:ssp}
Vladimir~G. Maz{'}ya, \emph{Sobolev spaces}, Springer Series in Soviet
  Mathematics, Springer-Verlag, Berlin, 1985, Translated from the Russian by T.
  O. Shaposhnikova. \MR{87g:46056}

\bibitem{mazya:97:dfb}
Vladimir~G. Maz{'}ya and Sergei~V. Poborchi, \emph{Differentiable functions on
  bad domains}, World Scientific Publishing Co. Inc., River Edge, NJ, 1997.
  \MR{99k:46057}

\bibitem{munkres:00:top}
James~R. Munkres, \emph{{Topology. 2nd ed.}}, {Upper Saddle River, NJ: Prentice
  Hall. xvi, 537 p. }, 2000.

\bibitem{oden:96:afa}
J.Tinsley Oden and Leszek~F. Demkowicz, \emph{{Applied functional analysis.}},
  {Boca Raton, FL: CRC Press. 653 p. }, 1996.

\bibitem{royden:88:ran}
H.L. Royden, \emph{{Real analysis. 3rd ed.}}, {New York: Macmillan Publishing
  Company; London: Collier Macmillan Publishing. xx, 444 p. }, 1988.

\bibitem{shvartsman:07:eos}
P.~Shvartsman, \emph{On extensions of {S}obolev functions defined on regular
  subsets of metric measure spaces}, J. Approx. Theory \textbf{144} (2007),
  no.~2, 139--161. \MR{MR2293385 (2007k:46057)}

\bibitem{swanson:99:sfi}
David Swanson and William~P. Ziemer, \emph{Sobolev functions whose inner trace
  at the boundary is zero}, Ark. Mat. \textbf{37} (1999), no.~2, 373--380.
  \MR{2000g:46048}

\bibitem{yosida:80:faa}
K{\^o}saku Yosida, \emph{Functional analysis}, 6th ed., Grundlehren der
  Mathematischen Wissenschaften, vol. 123, Springer-Verlag, Berlin, 1980.
  \MR{82i:46002}

\end{thebibliography}
\bibliographystyle{amsplain}

\end{document}